\documentclass[graybox]{svmult}

\usepackage{helvet}         
\usepackage{courier}        
\usepackage{type1cm}        
\usepackage{makeidx}         
\usepackage{graphicx}        
\usepackage{multicol}        
\usepackage[bottom]{footmisc}

\usepackage[toc,page]{appendix}

\usepackage{amsfonts}
\usepackage{amssymb}
\usepackage{amsmath}

\usepackage{cite}
\usepackage{epsf}

\def\beq{\begin{equation}}
\def\eqn#1{\beq\label{#1}}
\def\eeq{\end{equation}}

\def\eqna#1{\begin{eqnarray}\label{#1}}
\def\een{\end{eqnarray}}

\def\md{\medskip}
\def\nd{\end{document}}
\def\veps{\varepsilon}

\def\hh{h}

\def\nn{\nonumber}
\def\nt{\noindent}

\def\spa{\spadesuit}

\def\han{{\textstyle\frac{p+q-2}{2}}}

\def\ha{{\textstyle{\frac{1}{2}}}}

  \def\tV{{\tilde V}}
\def\tcn{{\tilde{\cal N}}}

\def\bu{\noindent $\bullet~$}

\def\deg{{\rm deg}\,}

\def\riga{-\kern-4pt - \kern-4pt -}
\font\fat=cmsy10 scaled\magstep5

\def\Bbullet{\raise-3pt\hbox{\fat\char"0F}}

\def\mt{\mapsto}

\def\down{\raise1.5pt\hbox{$\phantom{a}_2$}\downarrow}

\def\downa{\raise1.5pt\hbox{$\phantom{a}_{2\atop m_2}$}\downarrow}

\def\llr{\longrightarrow}

\def\({\left(}
\def\){\right)}
\def\eps{\epsilon}
 
\def\lra{\longrightarrow}

\def\dia{{$\diamondsuit$}}
 
\def\vf{\varphi}

\def\bbz{\mathbb{Z}}
\def\bbc{\mathbb{C}}
\def\bac{\bbc} 
\def\bbn{\mathbb{N}}
\def\a{\alpha}
\def\b{\beta}

\def\vr{\vert}

\def\l{\lambda}

\def\ca{{\cal A}}  \def\cc{{\cal C}}
\def\cd{{\cal D}}  \def\cf{{\cal F}}
\def\cg{{\cal G}} \def\ch{{\cal H}} 
  
\def\cm{{\cal M}} \def\cn{{\cal N}} 
\def\cp{{\cal P}}  
 \def\ct{{\cal T}}

\def\ido{intertwining differential operator}
\def\idos{intertwining differential operators}

\def\L{\Lambda}
\def\r{\rho}

\input epsf.tex
\newcount\figno
\def\fig#1#2#3{
\par\begingroup\parindent=0pt\leftskip=1cm\rightskip=1cm\parindent=0pt
\baselineskip=11pt \global\advance\figno by 1 
\epsfxsize=#3 \centerline{\epsfbox{#2}} \vskip 12pt
#1\par
\endgroup\par}
\def\figlabel#1{\xdef#1{\the\figno}}
\def\encadremath#1{\vbox{\hrule\hbox{\vrule\kern8pt\vbox{\kern8pt
\hbox{$\displaystyle #1$}\kern8pt} \kern8pt\vrule}\hrule}}

\begin{document}

\title*{Special Reduced Multiplets and Minimal
Representations for SO(p,q)}

\author{V.K. Dobrev}
\authorrunning{Dobrev}
\institute{V.K. Dobrev  \at
 Institute for Nuclear Research and Nuclear Energy,
 Tsarigradsko Chaussee 72, BG-1784 Sofia, Bulgaria,
 \email{dobrev@inrne.bas.bg}}

 \maketitle

\abstract{Using our previous results on the systematic construction of
invariant differential operators for non-compact semisimple Lie
groups we classify the special reduced multiplets and minimal
representations    in the case of SO(p,q).}


\section{Introduction}

In a recent paper \cite{Dobinv} we started the systematic explicit construction
of invariant differential operators. We gave an explicit description
of the building blocks, namely, the {\it parabolic subgroups and
subalgebras} from which the necessary representations are induced.
Thus we have set the stage for study of different non-compact
groups.

Since the study and description of detailed classification should be
done group by group we had to decide which groups to study first. We
decided to start with a subclass of   the hermitian-symmetric
algebras which share some special properties of the conformal
algebra ~$so(n,2)$. That is why, in view of applications to physics,
we called these algebras
  '{\it conformal Lie algebras}' (CLA), (or groups)
  \cite{Dobeseven}{}.
Later we gave a natural way to go beyond this subclass using essentially
the same results. For this we introduce the new notion of ~{\it
parabolic relation}~ between two non-compact semisimple Lie algebras
$\cg$ and $\cg'$ that have the same complexification and possess
maximal parabolic subalgebras with the same complexification
\cite{Dobparel}{}.

Thus, for example, using results for the conformal algebra
~$so(n,2)$~ (for fixed $n$) we can obtain results for all
pseudo-orthogonal algebras ~$so(p,q)$~ such that ~$p+q=n+2$. In this
way, in \cite{Dobparel} (among other things)  we gave the main and
the reduced multiplets of indecomposable elementary representations
for $so(p,q)$ including the necessary data for all relevant
invariant differential operators. We specially stressed that the
classification of all invariant differential operators includes as
special cases all possible ~{\it conservation laws}~ and ~{\it
conserved currents}, unitary or not. In the present paper we give
explicitly the conservation laws in the case of ~$so(p,q)$.

This paper is a short sequel of \cite{Dobparel}.
Due to the lack of space we refer to \cite{Dobparel}
for motivations and extensive list of literature on the subject.

\section{Preliminaries}

\nt Let ~$\cg=so(p,q)$, ~$p\geq q$, ~$p+q>4$.
 We choose a maximal parabolic $\cp = \cm \oplus \ca \oplus \cn$
 such that:
\begin{equation}\label{cmsosmay} \cm ~=~  so(p-1,q-1),\quad   \dim\, \ca ~=~1 \  ,\quad   \dim\, \cn ~=~   p+q
-2 \ .\end{equation} With this choice we get for the conformal algebra $so(n,2)$
the Bruhat decomposition ~$\cg = \cp \oplus \tcn$~ with direct  physical meaning
($\tcn\cong \cn$) \cite{Dobparel}{}.

 We label   the signature of the representations of
$\cg$   as follows: \begin {eqnarray}\label{sgnd}  &&\chi ~~=~~ \{\,
n_1,\, \ldots,\, n_{\hh}\ ; \  c\, \} \ , \\  &&\quad n_j \in \bbz/2\
, \quad c=d-\han\ , \quad \hh \equiv [\han] , \nn\\  && 0\leq \vr n_1 \vr
< n_2 < \cdots <  n_{\hh}\ , \quad p+q ~~{\rm even}\ , \nn\\  && 0 <
n_1 < n_2 < \cdots <  n_{\hh} \ , \quad p+q ~~{\rm odd}\ ,
\nn\end{eqnarray} where the parameter $c$ (related to the conformal weight $d$) labels the
characters of $\ca$,  and the first $\hh$ entries are labels of the
finite-dimensional (nonunitary for $q\neq 1$) irreps ~$\mu$~ of ~$\cm\,$.

Following \cite{DMPPT} we call the above induced
representations ~$\chi =$ Ind$^G_{\cp}(\mu\otimes\nu \otimes 1)$~
~{\it  elementary representations} (ERs) of $G=SO(p,q)$.
 Their spaces of functions are:  \begin {eqnarray}\label{func}
\cc_\chi ~&=&~ \{ \cf \in C^\infty(G,V_\mu) ~ \vr ~ \cf (gman) ~=~
e^{-\nu(H)} \cdot D^\mu(m^{-1})\, \cf (g) \} \nn\end{eqnarray}  where ~$a=
\exp(H)$, ~$H\in\ca\,$, ~$m\in M=SO(p-1,q-1)$, ~$n\in N=\exp\cn$. The
representation action is the {\it left regular action}:  \begin{equation}\label{lrega}
(\ct^\chi(g)\cf) (g') ~=~ \cf (g^{-1}g') ~, \quad g,g'\in G\ .\end{equation}
{\bf Remark:} ~Note that the group ~$M$~ has more general irreps
representing the centre of ~$M$. However, these are discrete parameters
which are not essential for the classification of the reducible ERs, cf.
\cite{Dobpq,Dob}.\dia

\md

\bu An important ingredient in our considerations are the ~{\it \it
highest/lowest weight representations}~ of ~$\cg^\bac$. These can be
realized as (factor-modules of) Verma modules ~$V^\L$~ over
~$\cg^\bac$, where ~$\L\in (\ch^\bac)^*$, ~$\ch^\bac$ is a Cartan
subalgebra of ~$\cg^\bac$, weight ~$\L = \L(\chi)$~ is determined
uniquely from $\chi$ \cite{Dob}{}.

Actually, since our ERs are induced from finite-dimensional
representations of ~$\cm$~   the Verma modules are
always reducible. Thus, it is more convenient to use ~{\it
generalized Verma modules} ~$\tV^\L$~ such that the role of the
highest/lowest weight vector $v_0$ is taken by the
(finite-dimensional) space ~$V_\mu\,v_0\,$. For the generalized
Verma modules (GVMs) the reducibility is controlled only by the
value of the conformal weight $d$, or the parameter $c$. Relatedly, for the \idos{} only
the reducibility w.r.t. non-compact roots is essential. Thus, from now on we shall consider
the ERs factored by the maximal invariant subspace generated by reducibilities w.r.t.
compact roots. We shall call these factored ERS: ~{\it compactly restricted ERs}.

\md

\bu One main ingredient of our approach is as follows. We group the
(reducible) ERs with the same Casimirs in sets called ~{\it
multiplets} \cite{Dobpq}{}. The multiplet corresponding to fixed
values of the Casimirs may be depicted as a connected graph, the
{\it vertices} of which correspond to the reducible ERs and the {\it
lines (arrows)}  between the vertices correspond to intertwining
operators. The multiplets contain explicitly all the data necessary
to construct the \idos{}. Actually, the data for each \ido{}
consists of the pair ~$(\b,m)$, where $\b$ is a (non-compact)
positive root of ~$\cg^\bac$, ~$m\in\bbn$, such that the {\it BGG
Verma module reducibility condition} \cite{BGG} (for highest weight
modules) is fulfilled:
\begin{equation}\label{bggr} (\L+\r, \b^\vee ) ~=~ m \ , \quad \b^\vee \equiv 2 \b
/(\b,\b) \ \end{equation} where $\r$ is half the sum of the positive roots of
~$\cg^\bac$. When the above holds then the Verma module with shifted
weight ~$V^{\L-m\b}$ (or ~$\tV^{\L-m\b}$ ~ for GVM and $\b$
non-compact) is embedded in the Verma module ~$V^{\L}$ (or
~$\tV^{\L}$). This embedding is realized by a singular vector
~$v_s$~  expressed by a polynomial ~$\cp_{m,\b}(\cg^-)$~ in the universal
enveloping algebra ~$(U(\cg_-))\ v_0\,$, ~$\cg^-$~ is the subalgebra
of ~$\cg^\bac$ generated by the negative root generators \cite{Dix}{}.
 More explicitly, \cite{Dob} ~$v^s_{m,\b} = \cp_{m,\b}\, v_0$ (or ~$v^s_{m,\b} =
 \cp_{m,\b}\, V_\mu\,v_0$ for GVMs).\footnote{For
explicit expressions for singular vectors we refer to \cite{Dobsin}{}.}\\
   Then there exists \cite{Dob} an ~{\it \ido{}}~ of order ~$m=m_\b$~:
  \eqn{invop}
\cd^m_{\b} ~:~ \cc_{\chi(\L)} ~\llr ~ \cc_{\chi(\L-m\b)} \eeq given
explicitly by: \eqn{singvv}\cd^m_{\b} ~=~ \cp^m_{\b}(\widehat{\cg^-})
\eeq where ~$\widehat{\cg^-}$~ denotes the $right$ action on the
functions ~$\cf$, cf. \eqref{fun}.

Thus, in each such situation we have an ~{\it invariant differential equation}~ of order ~$m=m_\b$~:
\begin{equation}\label{invde} \cd^m_{\b}\ f ~=~ f' \ , \qquad f \in \cc_{\chi(\L)} \ , \quad
f' \in \cc_{\chi(\L-m\b)} \ .\end{equation}

In many such situations the invariant operator ~$\cd^m_{\b}$~ has a non-trivial invariant
kernel.  These kernels are very important since in them
are realized the (irreducible) subrepresentations of $\cg$ as solutions of the equations:
\begin{equation}\label{invdec} \cd^m_{\b}\ f ~=~ 0 \ , \qquad f \in \cc_{\chi(\L)} \ ,
   \end{equation}
Furthermore, in some physical applications in the case of first order differential operators,
i.e., for ~$m=m_\b = 1$, equations \eqref{invdec}
are called ~{\it conservation laws}, and the elements ~$f\in \ker \cd_{m,\b}$~
are called ~{\it conserved currents}.

\section{Classification of  reducible ERs  for  ~{\boldmath $so(p,q)$}}

The reducible ERs are grouped in various multiplets. We start with the so-called
 main multiplets (which contain the maximal number of ERs
with this parabolic). We present them
with the following explicit parametrization of
the ERs in the multiplets (following \cite{Dobsrni}, see also \cite{DoPe:78}):
\begin{eqnarray}\label{sgne} \chi^\pm_1 &=& \{ \eps\, n_1,\, \ldots,\,
n_\hh \ ; \  \pm n_{\hh+1} \} \ ,
 \quad n_\hh < n_{\hh+1}\ , \nn\\
\chi^\pm_2 &=& \{ \eps\,  n_1,\, \ldots,\, n_{\hh-1},\,
n_{\hh+1}\ ; \  \pm n_\hh \}     \nn\\
\chi^\pm_3  &=&  \{  \eps\,
n_1,  \ldots,  n_{\hh-2},  n_{\hh},  n_{\hh+1}\ ; \  \pm
n_{\hh-1} \}     \nn\\  ... \\
 \chi^\pm_{\hh-1} &=& \{ \eps\, n_1,\,
n_2,\, n_4,\,\ldots,\, n_{\hh},\, n_{\hh+1}\ ; \  \pm  n_3 \}  \nn\\
 \chi^\pm_{\hh} &=& \{ \eps\, n_1,\,
n_3,\, \ldots,\, n_{\hh},\, n_{\hh+1}\ ; \  \pm  n_2 \}     \nn\\
\chi^\pm_{\hh+1} &=& \{ \eps\, n_2,\, n_3,\, \ldots,\, n_{\hh},\,
n_{\hh+1}\ ; \  \pm  n_1 \}     \nn\\  &\eps =& \begin{cases}
 \pm,\, ~&~  p+q ~~ even  \nn\\
                     1,  ~&~  p+q  ~~ odd \end{cases} \nn\end{eqnarray}
($\eps = \pm$~ is correlated with $\chi^\pm$).
Clearly, the multiplets correspond 1-to-1 to the finite-dimensional irreps
of $so(p+q,\bbc)$ with signature $\{n_1,\ldots,n_h,n_{h+1}\}$ and we are able to use
previous results due to the parabolic relation between the $so(p,q)$ algebras
for $p+q$~-fixed.
Note that the two representations in each pair ~$\chi^\pm$~ are called ~{\it shadow fields}.

Further, we denote by ~$\cc^\pm_i$~ the representation space with signature
~$\chi^\pm_i\,$.

The ERs in the multiplet are related by {\it intertwining integral and
differential operators}.

The  {\it integral operators} were introduced by
Knapp and Stein \cite{KnSt}{}.
Here these operators intertwine the pairs ~$\cc^\pm_i$~ (cf.
\eqref{sgne}): \begin{equation}\label{knapps} G^\pm_i ~:~ \cc^\mp_i \lra \cc^\pm_{i}
\ , \quad i ~=~ 1,\ldots,h+1  \ . \end{equation}

The {\it \idos}\ correspond to non-com\-pact positive roots of the
root system of ~$so(p+q,\bbc)$, cf. \cite{Dob}{}. In the current
context, compact roots of $so(p+q,\bbc)$ are those that are roots
also of the subalgebra $\cm^\bbc$, the rest of the roots are
non-compact. We denote the differential operators by
~$d_i,\,d'_i\,$. The spaces from \eqref{sgne} they intertwine   are:
\begin {eqnarray}\label{idoseven} && d_i ~:~ \cc^-_i ~\lra~ \cc^-_{i+1} \ , \quad i = 1,\ldots,h ~ ;\nn\\
&& d'_i ~:~ \cc^+_{i+1} ~\lra~ \cc^+_{i} \ , \quad i = 1,\ldots,h-1 ~ ;\nn\\
&& d_h ~:~ \cc^+_{h+1} \lra \cc^+_{h} \ , \quad (p+q)-{\rm even}~; \nn\\
&& d'_h ~:~ \cc^-_{h} \lra \cc^+_{h+1} \ , \quad (p+q)-{\rm even}~; \nn\\
&& d'_h ~:~ \cc^-_{h+1} \lra \cc^+_{h} \ , \quad (p+q)-{\rm even}~; \nn\\
&& d'_{h} ~:~ \cc^+_{h+1} \lra \cc^+_{h} \ , \quad (p+q)-{\rm odd}~; \nn\\
&& d_{h+1} ~:~ \cc^-_{h+1} \lra \cc^+_{h+1} \ , \quad (p+q)-{\rm odd}~ .
 \end{eqnarray}
The degrees of these \idos\ are given just by the differences of the
~$c$~ entries \cite{Dobsrni}{}: \begin {eqnarray}\label{degr} &&\deg
d_i = \deg d'_i = n_{\hh+2-i} - n_{\hh+1-i} ~=~ m_{h+2-i},\, \qquad i =
1,\ldots,h ,\,\nn\\  &&
\deg d'_{h} = n_2 + n_1 ~=~  m_1
,\, \quad (p+q) - {\rm even}\ , \nn\\  &&
\deg d_{h+1} = 2n_1 = m'_{h+1}  ~=~ m_1 ,\, \quad (p+q) - {\rm odd} \ .
\end{eqnarray}
where $d'_h$ is omitted from the first line for $(p+q)$ even.

\section{Multiplets and representations for ~$p+q$~ odd}

\subsection{Reduced multiplets for ~$p+q$~ odd}

In this Section we consider the case ~$p+q$~ odd, thus ~$h = \ha (p+q-3)$.
First we rewrite the main multiplets from \eqref{sgne} in the following parametrization:
\begin{eqnarray}\label{sgneod} \chi^\pm_1 &=& [  m_1,\, \ldots,\,
m_\hh 
\ ; \  \pm \ha (m_{1} + 2 m_{2,h+1} ) \ ]  \ ,  \\
\chi^\pm_2 &=& [   m_1,\, \ldots,\, m_{\hh-1},\,
m_{h,\hh+1} 
    \ ; \  \pm \ha (m_{1}+ 2m_{2,h}  )  \ ]      \nn\\
\chi^\pm_3  &=&  [
m_1,  \ldots,  m_{\hh-2},  m_{h-1,\hh},  m_{\hh+1}
\ ; \  \pm
\ha (m_{1} + 2m_{2,h-1} ) \ ]      \nn\\
... \nn\\
\chi^\pm_i  &=&  [
m_1,  \ldots,  m_{h-i+1},  m_{h-i+2,h-i+3},  m_{h+4-i}, \ldots,
m_{\hh},\, m_{\hh+1}\ ; \nn\\ &&\qquad   \pm
\ha (m_{1} + 2m_{2,h+2-i} ) \ ]      \nn\\
... \nn\\
 \chi^\pm_{\hh-1} &=& [  m_1,\,
m_2,\, m_{34},\,m_5,\,\ldots,\, m_{\hh},\, m_{\hh+1}
\ ; \
\pm \ha (m_{1} + 2m_{2,3} )  \ ]   \nn\\
 \chi^\pm_{\hh} &=& [  m_1,\,
m_{23},\, m_4,\,\ldots,\, m_{\hh},\, m_{\hh+1}
\ ; \  \pm  \ha (m_{1} + 2m_2 ) \ ]      \nn\\
\chi^\pm_{\hh+1} &=& [ m_1+ 2m_2,\, m_3,\, \ldots,\, m_{\hh},\,
m_{\hh+1}
\ ; \  \pm \ha m_1 \ ]      \nn\end{eqnarray}
where the last entry (as before) is the value of ~$c$,
while ~$m_i \in\bbn$~ are the Dynkin labels
(as in \eqref{degr}):
\eqn{dynkod}
  m_1 ~=~  2n_1 ~=~ 2 \ell_1+1   \ , \qquad
 m_{j} ~=~ n_{j} - n_{j-1} ~=~  \ell_{j} - \ell_{j-1} +1 ,\,  ~~j = 2,\ldots,h+1 \ .
\eeq
and we use the shorthand notation:
\eqn{shrt} m_{r,s} ~\equiv~ m_r + \cdots + m_s \ , ~~~ r<s \ ,
\qquad  m_{r,r} ~\equiv~ m_r \ ,
\qquad  m_{r,s} ~\equiv~ 0  \ , ~~~ r>s \ , \eeq
and we have also introduced the labels ~$\ell_k$~ (in order to facilitate comparison with the
literature):
\eqn{eljk} \ell_k ~=~ n_k - k +\ha \ , \qquad 0 \leq \ell_1 \leq \cdots \leq \ell_{h+1} \ . \eeq

We know that the ERs in a pair are related by the KS operators ~$G^\pm_i$~ \eqref{knapps},
however for ~$p+q$~ odd the operator ~$G^+_{h+1}$~ degenerates to a differential operator
of degree ~$m_1$~ corresponding to the only short
non-compact root ~$\veps_1\,$. The main multiplets are given   the Figure 1.
Note that following \cite{Dobparel}
we do not give the KS integral operators. Their presence is assumed by the symmetry w.r.t
the bullet in the centre of the Figure.

In this case there are ~$h+1$~ reduced multiplets which may be  obtained by
formally setting one Dynkin label to zero.
For ~$m_j=0$~ we denote the signatures by ~$^{j}\chi^\pm_k\,$.

We shall see that in every multiplet there is only one pair
(which we mark with $\spa$)
whose representations are of direct physical relevance (including finite-dimensional
irreps of the ~$\cm$~ subalgebra). Yet we list the others since they are related by
invariant differential operators which we record in each case.

In detail, the signatures are given similarly to \eqref{sgneod}:

\bu ~~$m_{h+1}=0$ ~equiv~ $n_{h+1} = n_h$
\begin{eqnarray}\label{sgneodr} ^{h+1}\chi^\pm_1 ~=~ ^{h+1}\chi^\pm_2 &=& [  m_1,\, \ldots,\,
m_\hh \ ; \  \pm \ha (m_{1} + 2 m_{2,h} ) \ ]  \ , \quad \spa \nn\\
^{h+1}\chi^\pm_3  &=&  [
m_1,  \ldots,  m_{\hh-2},  m_{h-1,\hh},  0\ ; \  \pm
\ha (m_{1} + 2m_{2,h-1} ) \ ]      \nn\\
... \\
^{h+1} \chi^\pm_{\hh-1} &=& [  m_1,\,
m_2,\, m_{34},\,m_5,\,\ldots,\, m_{\hh},\, 0 \ ; \
\pm \ha (m_{1} + 2m_{2,3} )  \ ]   \nn\\
^{h+1}\chi^\pm_{\hh} &=& [  m_1,\,
m_{23},\, m_4,\,\ldots,\, m_{\hh},\, 0\ ; \  \pm  \ha (m_{1} + 2m_2 ) \ ]      \nn\\
^{h+1}\chi^\pm_{\hh+1} &=& [ m_1+ 2m_2,\, m_3,\, \ldots,\, m_{\hh},\,
0 \ ; \  \pm \ha m_1 \ ]      \nn\end{eqnarray}
Here there are two differential operators involving physically relevant representations,
cf. Figure 2~:
\eqna{diffoneo} \cd^{m_h}_{\veps_1-\veps_{3}} ~&:&~
\cc^-_1 ~=~ \cc^-_2 ~\lra ~ \cc^-_3 \nn\\
 \cd^{m_h}_{\veps_1+\veps_{3}} ~&:&~
\cc^+_3  ~\lra ~  \cc^+_1 ~=~ \cc^+_2  \een

\md

\bu ~~$m_{h}=0$ ~equiv~ $n_{h} = n_{h-1}$
\begin{eqnarray}\label{sgneodrh} ^h\chi^\pm_1 &=& [  m_1,\, \ldots,\,
m_{\hh-1},\, 0 \ ; \  \pm \ha (m_{1} + 2 m_{2,h-1} + 2 m_{h+1}) \ ]  \ ,  \nn\\
^h\chi^\pm_2 ~=~ ^h\chi^\pm_3 &=& [   m_1,\, \ldots,\, m_{\hh-1},\,
m_{\hh+1}\ ; \  \pm \ha (m_{1}+ 2m_{2,h-1}  )  \ ]   , \quad \spa   \nn\\
... \\
 ^h\chi^\pm_{\hh-1} &=& [  m_1,\,
m_2,\, m_{34},\,m_5,\,\ldots,\, m_{\hh-1},\,0 ,\, m_{\hh+1}\ ; \
\pm \ha (m_{1} + 2m_{2,3} )  \ ]   \nn\\
 ^h\chi^\pm_{\hh} &=& [  m_1,\,
m_{23},\, m_4,\,\ldots,\, m_{\hh-1},\,0 ,\, m_{\hh+1}\ ; \  \pm  \ha (m_{1} + 2m_2 ) \ ]      \nn\\
^h\chi^\pm_{\hh+1} &=& [ m_1+ 2m_2,\, m_3,\, \ldots,\,m_{\hh-1},\, 0,\,
m_{\hh+1}\ ; \  \pm \ha m_1 \ ]      \nn\end{eqnarray}
Here there are four differential operators involving physically relevant representations,
cf. Figure 3~:
\eqna{difftwoo} \cd^{m_{h+1}}_{\veps_1-\veps_{2}} ~&:&~
\cc^-_1 ~\lra ~ \cc^-_2 ~=~ \cc^-_3 \nn\\
\cd^{m_{h-1}}_{\veps_1-\veps_{4}} ~&:&~
\cc^-_2 ~=~ \cc^-_3 ~\lra ~ \cc^-_4 \nn\\
 \cd^{m_{h-1}}_{\veps_1+\veps_{4}} ~&:&~
\cc^+_4  ~\lra ~  \cc^+_2 ~=~ \cc^+_3  \nn\\
\cd^{m_{h+1}}_{\veps_1+\veps_{2}} ~&:&~
 \cc^+_2 ~=~ \cc^+_3 \lra \cc^+_1 \een

The above case is typical for ~$m_k=0$~ for ~$k>2$. Then for ~$k= 2,1$~ we have:

\md



\bu ~~$m_{2}=0$ ~equiv~ $n_{2} = n_1$
\begin{eqnarray}\label{sgneodw} ^2\chi^\pm_1 &=& [  m_1,\,0,\,m_3,\, \ldots,\,
m_\hh \ ; \  \pm \ha (m_{1} + 2 m_{3,h+1} ) \ ]  \ ,  \nn\\
^2\chi^\pm_2 &=& [   m_1,\,0,\,m_3,\, \ldots,\, m_{\hh-1},\,
m_{h,\hh+1}\ ; \  \pm \ha (m_{1}+ 2m_{3,h}  )  \ ]      \nn\\
^2\chi^\pm_3  &=&  [
m_1,0,\,m_3,\,  \ldots,  m_{\hh-2},  m_{h-1,\hh},  m_{\hh+1}\ ; \  \pm
\ha (m_{1} + 2m_{3,h-1} ) \ ]      \nn\\
... \\
 ^2\chi^\pm_{\hh-1} &=& [  m_1,\,
0,\, m_{34},\,m_5,\,\ldots,\, m_{\hh},\, m_{\hh+1}\ ; \
\pm \ha (m_{1} + 2m_{3} )  \ ]   \nn\\
 ^2\chi^\pm_{\hh} ~=~ ^2\chi^\pm_{\hh+1} &=& [  m_1,\,
m_{3},\, m_4,\,\ldots,\, m_{\hh},\, m_{\hh+1}\ ; \  \pm  \ha m_{1} \ ]    , \quad \spa   \nn\end{eqnarray}
Here there are three differential operators involving physically relevant representations,
cf. Figure 4~:
\eqna{difffoo} \cd^{m_{3}}_{\veps_1-\veps_{h}} ~&:&~
\cc^-_{h-1} ~\lra ~ \cc^-_h ~=~ \cc^-_{h+1} \nn\\
\cd^{m_{1}}_{\veps_1} ~&:&~
\cc^-_h ~=~ \cc^-_{h+1} ~\lra ~  \cc^+_h ~=~ \cc^+_{h+1} \nn\\
\cd^{m_{3}}_{\veps_1+\veps_{h}} ~&:&~
\cc^+_h ~=~ \cc^+_{h+1} ~\lra ~ \cc^+_{h-1}
\een

\md

\bu ~~$m_{1}=0$ ~equiv~ $n_{1} = 0$
\begin{eqnarray}\label{sgneodo} ^1\chi^\pm_1 &=& [0,\,  m_2,\, \ldots,\,
m_\hh \ ; \  \pm   m_{2,h+1}  \ ]  \ ,  \nn\\
^1\chi^\pm_2 &=& [ 0,\,  m_2,\, \ldots,\, m_{\hh-1},\,
m_{h,\hh+1}\ ; \  \pm m_{2,h}    \ ]      \nn\\
^1\chi^\pm_3  &=&  [
0,\, m_2,  \ldots,  m_{\hh-2},  m_{h-1,\hh},  m_{\hh+1}\ ; \  \pm
 m_{2,h-1}  \ ]      \nn\\
 ... \\
 ^1\chi^\pm_{\hh-1} &=& [  0,\,
m_2,\, m_{34},\,m_5,\,\ldots,\, m_{\hh},\, m_{\hh+1}\ ; \
\pm m_{2,3}   \ ]   \nn\\
 ^1\chi^\pm_{\hh} &=& [  0,\,
m_{23},\, m_4,\,\ldots,\, m_{\hh},\, m_{\hh+1}\ ; \  \pm  m_2  \ ]      \nn\\
^1\chi_{\hh+1} &=& [ 2m_2,\, m_3,\, \ldots,\, m_{\hh},\,
m_{\hh+1}\ ; \  0 \ ]  , \quad \spa    \nn\end{eqnarray}
Here there are two differential operators involving physically relevant representations,
cf. Figure 5~:
\eqna{difffoz} \cd^{m_{2}}_{\veps_1-\veps_{h+1}} ~&:&~
\cc^-_{h-1} ~\lra ~ \cc^+_{h+1} ~=~ \cc^-_{h+1} \nn\\
 \cd^{m_{2}}_{\veps_1+\veps_{h+1}} ~&:&~
\cc^+_{h+1} ~=~ \cc^-_{h+1} ~\lra ~ \cc^+_{h-1} \ .
\een

For future reference we summarize the pairs of direct physical relevance reparametrizing
for more natural presentation and introducing uniform notation ~$_r\chi^\pm_k$~:
\eqn{sgnercol}
\begin{array}{ccccl}  _r\chi^\pm_1 ~&=&~ ^{h+1}\chi^\pm_1 &=& [  m_1,\, \ldots,\,
m_\hh \ ; \  \pm \ha (m_{1} + 2 m_{2,h} ) ] \ , \\
&&&&  d^+ \geq 2h   \ , ~~d^- \leq 1 \ , \\
_r\chi^\pm_2 ~&=&~ ^h\chi^\pm_2 &=& [   m_1,\, \ldots,\,
m_{\hh}\ ; \  \pm \ha (m_{1}+ 2m_{2,h-1}   )  ]  , \\
&&&&  d^+ \geq 2h-1    , ~~d^- \leq 2  ,  \\
... \\
_r\chi^\pm_j ~&=&~ ^{h-j+2}\chi^\pm_j  &=&  [
m_1,  \ldots,    m_{\hh}\ ; \  \pm
\ha (m_{1} + 2m_{2,h+1-j}  ) ]  \ , \\
&&&&  d^+ \geq 2h-j+1   \ , ~~d^- \leq j \ ,
~~~1\leq j \leq h-1\\
... \\
_r\chi^\pm_{h-1} ~&=&~ ^3\chi^\pm_{\hh-1} &=& [  m_1,\, \ldots,\, m_{\hh}\ ; \
\pm \ha (m_{1} + 2m_{2} )  ] \ , \\
&&&&  d^+ \geq h+2   \ , ~~d^- \leq h-1 \ , \\
_r\chi^\pm_h ~&=&~^2\chi^\pm_{\hh} &=& [  m_1,\, \ldots,\, m_{\hh}\ ; \  \pm  \ha m_{1}  ]  \ ,
\quad d^+ \geq h+1   \ , ~~d^- \leq h \ ,  \\
_r\chi_{h+1} ~&=&~ ^1\chi_{\hh+1} &=& [  2m_1,\, m_2,\, \ldots,\, m_{\hh}\ ; \  0 ]  \ ,
\quad d = h+\ha \
 \end{array}\eeq
where we have introduced notation ~$d^\pm$~ corresponding to the ~$"\pm"$~ occurrences:
\eqn{plumin} d^\pm ~=~ h + \ha  \pm |c| \ . \eeq

\subsection{Special reduced multiplets for ~$p+q$~ odd}

In addition to the standardly reduced multiplets discussed in the previous subsection,
there are special reduced multiplets   which may be formally obtained by formally setting one or two
Dynkin labels to a positive half integer. Again from each main multiplet only one pair is of physical relevance
but unlike the standardly reduced multiplets discussed in the previous subsection
these pairs are not related by differential operators
to the rest of the reduced multiplet (though having the same Casimirs). Thus, we present only the
physically relevant pairs.

\md

\bu ~~$m_{h+1}\mt \ha\mu\ , ~~~\mu \in 2\bbn-1$
\eqn{hone} _s\chi^\pm_1 = [  m_1,\, \ldots,\,
m_\hh
\ ; \  \pm \ha (m_{1} + 2 m_{2,h} + \mu ) \ ]
\eeq

\bu ~~$m_{h}\mt \ha\mu\ , ~~m_{h+1}\mt \ha\mu'\ , ~~~\mu,\mu' \in 2\bbn-1$
\eqn{htwo}
_s\chi^\pm_2 = [   m_1,\, \ldots,\, m_{\hh-1},\, \ha (\mu + \mu')
    \ ; \  \pm \ha (m_{1}+ 2m_{2,h-1} + \mu  )  \ ]      \eeq

\bu ~~$m_{h-1}\mt \ha\mu\ , ~~m_{h}\mt \ha\mu'\ , ~~~\mu,\mu' \in 2\bbn-1$
\eqn{hthree}
_s\chi^\pm_3  =  [m_1,  \ldots,  m_{\hh-2}, \ha (\mu + \mu') ,  m_{\hh+1}
\ ; \  \pm
\ha (m_{1} + 2m_{2,h-2} + \mu ) \ ]      \eeq

\bu ~~$m_{h-j+2}\mt \ha\mu\ , ~~m_{h-j+3}\mt \ha\mu'\ , ~~~\mu,\mu' \in 2\bbn-1\ ,
\quad 2\leq j \leq h$
\eqna{hijk}
_s\chi^\pm_j  &=&  [
m_1,  \ldots,  m_{h-j+1},  \ha (\mu + \mu')  ,  m_{h+4-j}, \ldots,
m_{\hh},\, m_{\hh+1}\ ; \ \nn\\ && \pm
\ha (m_{1} + 2m_{2,h+1-j} + \mu ) \ ] \
\een

\bu ~~$m_{3}\mt \ha\mu\ , ~~m_{4}\mt \ha\mu'\ , ~~~\mu,\mu' \in 2\bbn-1$
\eqn{hminusone}
_s\chi^\pm_{\hh-1} = [  m_1,\,
m_2,\,  \ha (\mu + \mu')  ,\,m_5,\,\ldots,\, m_{\hh},\, m_{\hh+1}
\ ; \ \pm \ha (m_{1} + 2m_{2} + \mu )  \ ]\eeq

\bu ~~$m_{2}\mt \ha\mu\ , ~~m_{3}\mt \ha\mu'\ , ~~~\mu,\mu' \in 2\bbn-1$
\eqn{hhhh}
_s\chi^\pm_{\hh} = [  m_1,\,  \ha (\mu + \mu') ,\, m_4,\,\ldots,\, m_{\hh},\, m_{\hh+1}
\ ; \  \pm  \ha (m_{1} + \mu ) \ ]      \eeq

\bu ~~$m_{2}\mt \ha\mu\ , ~~~\mu \in 2\bbn-1$
\eqn{hplusone}
_s\chi^\pm_{\hh+1} = [ m_1+ \mu ,\, m_3,\, \ldots,\, m_{\hh},\,
m_{\hh+1}\ ; \  \pm \ha m_1 \ ]      \eeq

In each pair there are the standard KS integral operators ~$G^\pm_k$~ between ~$_r\chi^\mp_k\,$.~
Furthermore, the ERs in a pair are reducible w.r.t. the compact roots
and in addition the ERs    ~$_r\chi^-_k$~ are reducible w.r.t. the
only short noncompact root ~$\veps_1\,$. Actually, the corresponding differential operators
are degenerations of the corresponding KS operators
  ~$G^+_k$~ \eqref{knapps}.
  (In the main multiplets the same   was happening but only for ~$k=h+1$.)
  Thus, we have:
\eqn{diffodd} \cd^{2|c_k|}_{\veps_1} ~:~
_r\cc^-_k ~\lra~ _r\cc^+_k \ , \qquad G^+_k  ~\sim ~\cd^{2|c_k|}_{\veps_1} \eeq
where ~$c_k$~ is the value of ~$c$~ of the   ER ~$_r\chi^-_k\,$.

\bigskip

Finally, we give a doubly reduced case originating from \eqref{hplusone} setting ~$m_1=0$~:
\eqn{hplusonez}
_{rs}\chi^\pm_{\hh+1} = [  \mu ,\, m_2,\, \ldots,\, m_{\hh}\ ; \  0 \ ] \ ,
\quad m_k\in\bbn, \quad  \mu \in 2\bbn-1 \ .   \eeq
This is a singlet and the ER is reducible only w.r.t. the compact roots,
there are no non-trivial differential operators, thus,
the corresponding generalized Verma module and the compactly
restricted ER are irreducible.

\bigskip

\subsection{Special cases for ~$p+q$~ odd}

The ERS  ~$\chi^-_1$~ are the only ones in the multiplet that contain as irreducible subrepresentations
the finite-dimensional   irreducible representations of ~$\cg$. More precisely, the ER ~$\chi^-_1$~
contains the finite-dimensional   irreducible representation of ~$\cg$~
with signature ~$(m_1, \ldots, m_{h+1})$. (Certainly, the latter is non-unitary
except the case of the trivial one-dimensional  obtained for ~$m_i=1$, ~$\forall\, i$.)

Another important case is the ER with signature
~$\chi^+_1\,$. It contains a unitary discrete series representation  of $so(p,q)$ realized on an
invariant subspace ~$\cd$~ of
the ER ~$\chi^+_1\,$. That subspace is annihilated by the KS operator
~$G^-_1\,$,\ and is the image of the KS operator ~$G^+_1\,$.

Furthermore when ~$p>q=2$~ the invariant subspace ~$\cd$~ is the direct sum of two subspaces
~$\cd ~=~ \cd^+ \oplus \cd^-$, in which are realized a
 {\it holomorphic discrete series representation} and its conjugate
  {\it   anti-holomorphic discrete
series representation}, resp.
Note that the corresponding lowest weight GVM is infinitesimally
equivalent only to the holomorphic discrete series, while the
conjugate highest weight GVM is infinitesimally equivalent to the
anti-holomorphic discrete series.

Thus, the signatures of the (holomorphic) discrete series are:
\eqn{holo}
\chi^+_1 = [  m_1,\, \ldots,\, m_\hh
\ ; \  d ~=~ h +  \ha (m_{1}  +1) + m_{2,h} ~+~ \nu \ ]\ , \quad \nu\in\bbn
\eeq

More (non-holomorphic) discrete series representations are contained in
 ~$\chi^+_k$~ for ~$1<k\leq h$.

The next important case are the ~{\it limits of (holomorphic) discrete series}~
which are contained in the reduced case \eqref{sgneodr}:
\eqn{limho}
_r\chi^+_1 ~=~  [  m_1,\, \ldots,\,
m_\hh \ ; \ d ~=~ h +  \ha (m_{1}  +1) + m_{2,h} \ ] \eeq
(with conformal weight obtained from \eqref{holo} as "limit" for $\nu=0$).

Finally, we mention the so called ~{\it first reduction points} (FRP). For ~$q=2$~ these are the boundary values of ~$d$~ from below
of the positive energy UIRs.  Most of the FRPs  are
contained in ~$\chi^+_{\hh+1}$, cf. \eqref{sgneod}, which we give with suitable reparametrization:
\eqn{frp} \chi^+_{\hh+1} ~=~ [ m_1,\, m_2,\, \ldots,\, m_{\hh}\ ; \  d ~=~ h  +   \ha m_1 -\ha \ ]\ ,
\qquad m_1 \geq 3 \ .\eeq
The  FRP cases for ~$m_1=1,2$~ (with the same values of $d$ by specializing $m_1$)
are found in \eqref{sgneodw}, \eqref{sgneodo}, resp:
\begin{eqnarray}\label{holoone}
_r\chi^-_{\hh} ~&=&~ [  1,\, m_{2},\, \ldots,\, m_{\hh} \ ; ~  d ~=~ h  \ ]\ , \\
_r\chi_{\hh+1} ~&=&~ [ 2,\, m_2,\, \ldots,\, m_{\hh}\ ; ~ d ~=~ h +\ha   \ ] \ .
\end{eqnarray}

Finally, we give some discrete unitary points below the FRP which are found in
the special reduced ERs \eqref{hhhh} (used for ~$m_1=\mu=1$, ~$\mu' = 2m_2-1$),
and then \eqref{hplusone}
used first for ~$m_1=2$, ~$\mu = 2k-1$,
and then for ~$m_1=1$, ~$\mu = 2k-1$~:
\begin{eqnarray}\label{single}
_s\chi^-_{\hh} ~&=&~ [  1,\,  m_2, \,\ldots,\,  m_{\hh}
\ ; \  d ~=~ h - \ha   \ ] \ ,      \\
_s\chi^-_{\hh+1} ~&=&~ [ 2k+1  ,\, m_2,\, \ldots,\, m_{\hh}\
  ; \ d ~=~ h -\ha \ ] \ , \quad k\in\bbn \\
_s\chi^-_{\hh+1} ~&=&~ [ 2k  ,\, m_2,\, \ldots,\, m_{\hh}\
  ; \ d ~=~ h  \ ] \ , \quad k\in\bbn
  \end{eqnarray}

\subsection{Minimal irreps for ~$p+q$~ odd}

First we give the minimal irreps occurring in
standardly reduced multiplets displaying together only the
physically relevant representations:
\begin{eqnarray}\label{sgnmin}
_r\chi^-_1  &=& [
 \ 1,\, \ldots,\, 1\ ; ~   d ~=~ 1  \ ]\ , \\
&& _rL^-_1 ~=~ \{\ \vf \in\ _r\cc^-_1  ~:~
\cd^{1}_{\veps_1-\veps_{3}}\, \vf ~=~ 0 \ ,
\quad G^+_1\, \vf ~=~ 0 \ \} \ , \nn\\
%
%
... \nn\\
_r\chi^-_{j}  &=& [ \ 1,\, \ldots,\,1    \ ; \
 d ~=~ j    \ ] \ , \nn\\
&& _rL^-_j ~=~ \{\ \vf \in\ _r\cc^-_j  ~:~
\cd^{1}_{\veps_1-\veps_{j+2}}\, \vf ~=~ 0 \ ,
\quad G^+_j\, \vf ~=~ 0 \ \} \ , \nn\\
&&\qquad\qquad\qquad 1\leq j \leq h-1\ ,\nn\\
... \nn\\
  _r\chi^-_{\hh} &=& [ \ 1,\, \ldots,\, 1 \ ; \  d_{\rm FRP} ~=~ h  \ ]  \ ,    \nn\\
&&_rL^-_h ~=~ \{\ \vf \in\ _r\cc^-_h  ~:~
\cd^{1}_{\veps_1}\, \vf ~=~ 0 \  \ \} , \qquad
 G^+_h ~=~ \cd^{1}_{\veps_1}  \ , \nn\\
_r\chi_{\hh+1} &=& [\  2,\, 1,\, \ldots,\, 1 \ ; \  d_{\rm FRP} ~=~ h+\ha\ ]  \ ,
   \nn\\
&&_rL_{h+1} ~=~ \{\ \vf \in\ _r\cc^-_{h+1}  ~:~
\cd^{1}_{\veps_1+\veps_{h+1}}\, \vf ~=~ 0 \  \ \} \nn
\end{eqnarray}
(In the last case there is no KS operator since ~$c=0$.)

We see that for ~$h\geq 2$~ there are discrete unitary points ~{\it below}~ the FRPs.
 For fixed ~$h\geq 2$~ these are
in ~$_r\chi^-_j$~ with conformal weight ~$d=j$~ (and trivial $\cm$ inducing irreps)
for ~$j = 1,\ldots, h-1$. Furthermore, as evident from \eqref{sgnevsum}
for ~$h\geq 3$~  there are discrete unitary points below those displayed.
For fixed ~$h\geq 3$~ these are
in ~$_r\chi^-_j$~ with conformal weight ~$1\leq d <j$~ (and non-trivial ~$\cm$~ inducing irreps)
for ~$j = 2,\ldots, h-2$. It seems that all this picture is consistent with \cite{EHW}.
More details  will be given elsewhere.

Next we give the case of special reduced multiplets displaying together the
physically relevant representations:
\begin{eqnarray}\label{sgnminsp} _s\chi^-_1 &=& [\   1,\, \ldots,\,
1  \ ; \  d ~=~ \ha \ ] \ ,   \\
&&_sL^-_1 ~=~ \{\ \vf \in\ _s\cc^-_1  ~:~
\cd^{2h}_{\veps_1}\, \vf ~=~ 0 \  \ \} , \qquad
 G^+_1 ~=~ \cd^{2h}_{\veps_1}  \ , \nn\\
... \nn\\
_s\chi^-_j  &=&  [\  1,\, \ldots,\, 1 \  ; \
 d ~=~ j-\ha   \ ]  \ , \qquad 1\leq j \leq h  \nn\\
&&_sL^-_j ~=~ \{\ \vf \in\ _s\cc^-_j  ~:~
\cd^{2(h+1-j)}_{\veps_1}\, \vf ~=~ 0 \  \ \} , \qquad
 G^+_j ~=~ \cd^{2(h+1-j)}_{\veps_1}  \ , \nn\\
... \nn\\
_s\chi^-_h  &=&  [\  1,\, \ldots,\, 1 \  ; \
 d ~=~ h-\ha   \ ]  \ ,  \nn\\
&&_sL^-_h ~=~ \{\ \vf \in\ _s\cc^-_h  ~:~
\cd^{2}_{\veps_1}\, \vf ~=~ 0 \  \ \} , \qquad
 G^+_h ~=~ \cd^{2}_{\veps_1}  \ , \nn\\
_s\chi_{\hh+1} &=& [\   2,\, 1,\, \ldots,\, 1 \ ; \  d ~=~ h  \ ]  \ ,
     \nn\\
&&_sL^-_{h+1} ~=~ \{\ \vf \in\ _s\cc^-_{h+1}  ~:~
\cd^{1}_{\veps_1}\, \vf ~=~ 0 \  \ \} , \qquad
 G^+_{h+1} ~=~ \cd^{1}_{\veps_1}  \ , \nn
     \end{eqnarray}
Here all irreps are below the FRP.
The "most" minimal representations are the last two cases of \eqref{sgnminsp}.
For ~$h=1$, i.e., ~$so(3,2)$~ these are the so-called ~${\it singletons}$~
discovered by Dirac \cite{Dirac}.


\subsection{Singular vectors needed for the invariant differential operators}

The mostly used case is ~$\veps_1 ~=~ \a_1 + \cdots + \a_\ell\,$, ~$\ell = h+1$.
The corresponding singular vector of weight ~$m\veps_1$~ is given in (13)
\cite{Dobsin} (noting that this is an ~$sl(n)$~ formula
in quantum group setting, thus, one should take $q=1$):
\eqna{singone}
v^m_{\veps_1}
~&=&~ \sum_{k_1=0}^m \cdots \sum_{k_{\ell-1}=0}^m a_{k_1\dots
k_{\ell-1}} (X^-_1)^{m-k_1} \cdots
(X^-_{\ell -1})^{m-k_{\ell-1}}
\times \nn\\ && \times
~(X^-_\ell)^m  ~(X^-_{\ell-1})^{k_{\ell-1}}~
  \cdots
~(X^-_1)^{k_1} ~\otimes v_0 \ , \\
 a_{k_1\dots k_{\ell-1}} ~&=&~
(-1)^{k_1 + \cdots + k_{\ell-1}} ~a^\ell ~\left({m \atop
k_1}\right) \cdots ~\left({m \atop k_{\ell-1}}\right) ~\times
\nn\\ &&\times ~{[(\l + \r)(H^1)]\over [(\l + \r)(H^1) - k_1]}
\cdots {[(\l + \r)(H^{\ell-1})]\over [(\l + \r)(H^{\ell-1}) - k_{\ell-1}]} ,  \nn\een
where ~$X^\pm_k$~ are the simple root vectors, ~$H_k$~ are the long Chevalley Cartan elements
~$H_k = [X^+_k,X^-_k]$, $k<\ell$, ~$H^s ~=~  H_1 + H_2 + \cdots + H_s\,$, ~$\r$~ is the half-sum
of the positive roots.

Other cases are: ~$\veps_1 - \veps_j ~=~ \a_1 + \cdots + \a_{j-1}\,$.
 Clearly, one uses again
formula \eqref{singone} replacing ~$\ell \mt j-1$.

The last case is:
~$\veps_1 + \veps_\ell ~=~ \a_1 + \cdots + \a_{\ell-1} + 2\a_\ell\,$,
~$\ell = h+1$. The singular vector  is given in (19) of \cite{Dobsin}:
\eqna{singtwo}
&& v^m_{\veps_1 + \veps_\ell} ~=~ \sum_{k_1=0}^m \cdots
\sum_{k_{\ell-2}=0}^m \sum_{k_{\ell-1}=0}^{2m} ~b_{k_1 \dots
k_{\ell-1}} ~(X^-_1)^{m-k_1} \cdots
~(X^-_{\ell-2})^{m-k_{\ell-2}}\ \times \nn\\ && \times\ (X^-_\ell)^{2m-k_{\ell-1}}
(X^-_{\ell-1})^m ~(X^-_\ell)^{k_{\ell-1}}
~(X^-_{\ell-2})^{k_{\ell-2}} \cdots (X^-_1)^{k_1} ~\otimes v_0
~, \\  && b_{k_1 \dots k_{\ell-1}} ~=~ (-1)^{k_1 +
\cdots + k_{\ell-1}} ~b^\ell ~\left({m \atop k_1}\right)
\cdots ~\left({m \atop k_{\ell-2}}\right) ~\left({2m \atop
k_{\ell-1}}\right) \ \times \nn\\ &&\times\ {[(\l + \r)(H^1)]
\over [(\l + \r)(H^1) - k_1]} \dots ~{[(\l +
\r)(H^{\ell-2})] \over [(\l + \r)(H^{\ell-2}) -
k_{\ell-2}]} ~{[(\l + \r)(H_\ell)] \over [(\l +
\r)(H_\ell) - k_{\ell-1}]}  \nn \een

\section{Multiplets and representations for ~$p+q$~ even}

\subsection{Reduced multiplets for ~$p+q$~ even}

In this Section we consider the case ~$p+q$~ odd, thus ~$h = \ha (p+q-2)$.
First we introduce the Dynkin labels parametrization of the multiplets:
\begin{eqnarray}\label{sgnev} \chi^\pm_1 &=& [
 (m_1,\, \ldots,\, m_\hh)^\pm
\ ; \  \pm (\ha m_{12} +  m_{3,h+1} ) \ ]  \ ,  \\
\chi^\pm_2 &=& [   (m_1,\, \ldots,\, m_{\hh-1},\,
m_{h,\hh+1})^\pm
    \ ; \  \pm (\ha m_{12}+ m_{3,h}  )  \ ]      \nn\\
\chi^\pm_3  &=&  [
(m_1,  \ldots,  m_{\hh-2},  m_{h-1,\hh},  m_{\hh+1})^\pm
\ ; \  \pm
(\ha m_{12} + m_{3,h-1} ) \ ]      \nn\\
... \nn\\
\chi^\pm_j  &=&  [
(m_1,  \ldots,  m_{h-j+1},  m_{h-j+2,h-j+3},  m_{h+4-j}, \ldots,
m_{\hh},\, m_{\hh+1})^\pm \ ; \ \nn\\
&& \pm
(\ha m_{12} + m_{3,h+2-j} ) \ ]\ , \quad 2\leq j \leq h-1 \ ,       \nn\\
... \nn\\
 \chi^\pm_{\hh-1} &=& [  (m_1,\,
m_2,\, m_{34},\,m_5,\,\ldots,\, m_{\hh},\, m_{\hh+1})^\pm
\ ; \
\pm  (\ha m_{12} + m_{3} )  \ ]   \nn\\
 \chi^\pm_{\hh} &=& [  (m_{1'3},\,
m_{23},\, m_4,\,\ldots,\, m_{\hh},\, m_{\hh+1})^\pm
\ ; \  \pm  \ha m_{12}  \ ]      \nn\\
\chi^\pm_{\hh+1} &=& [ (m_{13},\, m_3,\, \ldots,\, m_{\hh},\,
m_{\hh+1})^\pm \ ; \  \pm \ha (m_1-m_2) \ ]      \nn\end{eqnarray}
where the conjugation of the $\cm$ labels interchanges the first two
entries: \eqna{conju} &&(m_1,\, \ldots,\, m_\hh)^- ~=~ (m_1,\, \ldots,\,
m_\hh), \\ && (m_1,\,m_2,\, \ldots,\, m_\hh)^+ ~=~ (m_2,\,m_1,\,
\ldots,\, m_\hh)\ , \nn\een the last entry (as before) is the value of
~$c$, while ~$m_i \in\bbn$~ are the Dynkin labels (as in
\eqref{degr}): \eqna{dynkod}
&&  m_1 ~=~  n_1 + n_2 ~=~  \ell_1+\ell_2 +1   \ , \\
&& m_{j} ~=~ n_{j} - n_{j-1} ~=~  \ell_{j} - \ell_{j-1} +1 ,\,  ~~j = 2,\ldots,h+1 \ ,
 \nn\een
finally, ~$m_{1'3} \equiv m_1+m_3$.

The main  multiplets are given in Figure 6. Note that as in the odd case
we do not give the KS integral operators.

\bigskip

Then we give the reduced multiplets:

\bu ~~$m_{h+1}=0$ ~equiv~ $n_{h+1} = n_h$
\begin{eqnarray}\label{sgnevhplus} \chi^\pm_1 ~=~ \chi^\pm_2 &=& [
 (m_1,\, \ldots,\, m_\hh)^\pm
\ ; \  \pm (\ha m_{12} +  m_{3,h} ) \ ]  \ ,\qquad  \spa \\
 \chi^\pm_3  &=&  [
(m_1,  \ldots,  m_{\hh-2},  m_{h-1,\hh},  0)^\pm
\ ; \  \pm
(\ha m_{12} + m_{3,h-1} ) \ ]      \nn\\
... \nn\\
\chi^\pm_i  &=&  [
(m_1,  \ldots,  m_{h-i+1},  m_{h-i+2,h-i+3},  m_{h+4-i}, \ldots,
m_{\hh},\, 0)^\pm \ ; \nn\\ &&\qquad \qquad  \pm
(\ha m_{12} + m_{3,h+2-i} ) \ ]      \nn\\
... \nn\\
 \chi^\pm_{\hh-1} &=& [  (m_1,\,
m_2,\, m_{34},\,m_5,\,\ldots,\, m_{\hh},\, 0)^\pm
\ ; \
\pm  (\ha m_{12} + m_{3} )  \ ]   \nn\\
 \chi^\pm_{\hh} &=& [  (m_{1'3},\,
m_{23},\, m_4,\,\ldots,\, m_{\hh},\, 0)^\pm
\ ; \  \pm  \ha m_{12}  \ ]      \nn\\
\chi^\pm_{\hh+1} &=& [ (m_{13},\, m_3,\, \ldots,\, m_{\hh},\,0)^\pm
\ ; \  \pm \ha (m_1-m_2) \ ]      \nn\end{eqnarray}
Here there are two differential operators involving physically relevant representations,
cf. Figure 7~:
\eqna{diffone} \cd^{m_h}_{\veps_1-\veps_{3}} ~&:&~
\cc^-_1 ~=~ \cc^-_2 ~\lra ~ \cc^-_3 \nn\\
 \cd^{m_h}_{\veps_1+\veps_{3}} ~&:&~
\cc^+_3  ~\lra ~  \cc^+_1 ~=~ \cc^+_2  \een

\md

\bu ~~$m_{h}=0$ ~equiv~ $n_{h} = n_{h-1}$
\begin{eqnarray}\label{sgnevh} \chi^\pm_1 &=& [
 (m_1,\, \ldots,\, m_{h-1},\, 0)^\pm
\ ; \  \pm (\ha m_{12} +  m_{3,h-1} + 2m_{h+1} ) \ ]  \ ,  \\
\chi^\pm_2 = \chi^\pm_3 &=& [   (m_1,\, \ldots,\, m_{\hh-1},\,
m_{\hh+1})^\pm
    \ ; \  \pm (\ha m_{12}+ m_{3,h-1}  )  \ ]  \, \qquad \spa    \nn\\
 ... \nn\\
\chi^\pm_i  &=&  [
(m_1,  \ldots,  m_{h-i+1},  m_{h-i+2,h-i+3},  m_{h+4-i}, \ldots,
m_{\hh-1},\, 0,\, m_{\hh+1})^\pm \ ; \nn\\ &&\qquad\qquad  \pm
(\ha m_{12} + m_{3,h+2-i} ) \ ]      \nn\\
... \nn\\
 \chi^\pm_{\hh-1} &=& [  (m_1,\,
m_2,\, m_{34},\,m_5,\,\ldots,\, m_{\hh-1},\, 0,\, m_{\hh+1})^\pm
\ ; \
\pm  (\ha m_{12} + m_{3} )  \ ]   \nn\\
 \chi^\pm_{\hh} &=& [  (m_{1'3},\,
m_{23},\, m_4,\,\ldots,\, m_{\hh-1},\,0,\, m_{\hh+1})^\pm
\ ; \  \pm  \ha m_{12} \ ]      \nn\\
\chi^\pm_{\hh+1} &=& [ (m_{13},\, m_3,\, \ldots,\, m_{\hh-1},\,0,\,
m_{\hh+1})^\pm
\ ; \  \pm \ha (m_1-m_2) \ ]      \nn\end{eqnarray}
Here there are four differential operators involving physically relevant representations,
cf. Figure 8~:
\eqna{difftwo} \cd^{m_{h+1}}_{\veps_1-\veps_{2}} ~&:&~
\cc^-_1 ~\lra ~ \cc^-_2 ~=~ \cc^-_3 \nn\\
\cd^{m_{h-1}}_{\veps_1-\veps_{4}} ~&:&~
\cc^-_2 ~=~ \cc^-_3 ~\lra ~ \cc^-_4 \nn\\
 \cd^{m_{h-1}}_{\veps_1+\veps_{4}} ~&:&~
\cc^+_4  ~\lra ~  \cc^+_2 ~=~ \cc^+_3  \nn\\
\cd^{m_{h+1}}_{\veps_1+\veps_{2}} ~&:&~
 \cc^+_2 ~=~ \cc^+_3 \lra \cc^+_1 \een

The above case is typical for ~$m_k=0$~ for ~$k>3$. Then for ~$k= 3,2,1$~ we have:

\bu ~~$m_{3}=0$ ~equiv~ $n_{3} = n_{2}$
\begin{eqnarray}\label{sgnevthree} \chi^\pm_1 &=& [
 (m_1,\,m_2,\,0,\, m_4,\, \ldots,\, m_\hh)^\pm
\ ; \  \pm (\ha m_{12} +  m_{4,h+1} ) \ ]  \ ,  \\
\chi^\pm_2 &=& [   (m_1,\, m_2,\,0,\, m_4,\,\ldots,\, m_{\hh-1},\,
m_{h,\hh+1})^\pm
    \ ; \  \pm (\ha m_{12}+ m_{4,h}  )  \ ]      \nn\\
... \nn\\
\chi^\pm_i  &=&  [
(m_1,\,m_2,\,0,\, m_4,\,  \ldots,  m_{h-i+1},  m_{h-i+2,h-i+3},  m_{h+4-i}, \ldots,
\nn\\ &&\qquad\qquad
m_{\hh},\, m_{\hh+1})^\pm \ ;     \pm
(\ha m_{12} + m_{4,h+2-i} ) \ ]      \nn\\
... \nn\\
 \chi^\pm_{\hh-1} = \chi^\pm_{\hh} &=& [  (m_1,\,
m_2,\, m_{4},\,\ldots,\, m_{\hh},\, m_{\hh+1})^\pm
\ ; \
\pm  \ha m_{12}   \ ]  \ ,\qquad \spa \nn\\
  \chi^\pm_{\hh+1} &=& [ (m_{12},0,\, \, m_4,\, \ldots,\, m_{\hh},\,
m_{\hh+1})^\pm
\ ; \  \pm \ha (m_1-m_2) \ ]      \nn\end{eqnarray}
Here there are six differential operators involving physically relevant representations,
cf. Figure 9~:
\eqna{diffthree} \cd^{m_{4}}_{\veps_1-\veps_{h-1}} ~&:&~
\cc^-_{h-2} ~\lra ~ \cc^-_h ~=~ \cc^-_{h-1} \nn\\
\cd^{m_{2}}_{\veps_1-\veps_{h+1}} ~&:&~
\cc^-_h ~=~ \cc^-_{h-1} ~\lra ~ \cc^-_{h+1} \nn\\
\cd^{m_{1}}_{\veps_1+\veps_{h+1}} ~&:&~
\cc^-_h ~=~ \cc^-_{h-1} ~\lra ~ \cc^+_{h+1} \nn\\
 \cd^{m_{1}}_{\veps_1+\veps_{h+1}} ~&:&~
\cc^-_{h+1}  ~\lra ~  \cc^+_h ~=~ \cc^+_{h-1}  \nn\\
\cd^{m_{2}}_{\veps_1-\veps_{h+1}} ~&:&~
 \cc^+_{h+1} ~\lra~  \cc^+_h ~=~ \cc^+_{h-1} \nn\\
\cd^{m_{4}}_{\veps_1+\veps_{h-1}} ~&:&~
\cc^+_h ~=~ \cc^+_{h-1} ~\lra~ \cc^+_{h-2}
  \een

\md

\bu ~~$m_{2}=0$ ~equiv~ $n_{2} = n_{1}$,
\begin{eqnarray}\label{sgnevtwo} _2\chi^\pm_1 &=& [
 (m_1,\,0,m_3,\, \ldots,\, m_\hh)^\pm
\ ; \  \pm (\ha m_1 +  m_{3,h+1} ) \ ]  \ ,  \\
_2\chi^\pm_2 &=& [   (m_1,\,0,m_3,\, \ldots,\, m_{\hh-1},\,
m_{h,\hh+1})^\pm
    \ ; \  \pm (\ha m_1+ m_{3,h}  )  \ ]      \nn\\
_2\chi^\pm_3  &=&  [
(m_1,0,m_3,\,  \ldots,  m_{\hh-2},  m_{h-1,\hh},  m_{\hh+1})^\pm
\ ; \  \pm
(\ha m_1 + m_{3,h-1} ) \ ]      \nn\\
... \nn\\
_2\chi^\pm_i  &=&  [
(m_1, 0,m_3,\, \ldots,  m_{h-i+1},  m_{h-i+2,h-i+3},  m_{h+4-i}, \ldots,
\nn\\ &&\qquad\qquad
m_{\hh},\, m_{\hh+1})^\pm \ ; \  \pm
(\ha m_1 + m_{3,h+2-i} ) \ ]      \nn\\
... \nn\\
 _2\chi^\pm_{\hh-1} &=& [  (m_1,\,
0,\, m_{34},\,m_5,\,\ldots,\, m_{\hh},\, m_{\hh+1})^\pm
\ ; \
\pm  (\ha m_1 + m_{3} )  \ ]   \nn\\
 _2\chi^\pm_{\hh} = _2\chi^\pm_{\hh+1} &=& [  (m_1+m_{3},\,
m_{3},\,\ldots,\, m_{\hh},\, m_{\hh+1})^\pm
\ ; \  \pm  \ha m_1  \ ] \ , \qquad \spa
       \nn\end{eqnarray}
Here there are three differential operators involving physically relevant representations,
cf. Figure 10~:
\eqna{difffo} \cd^{m_{3}}_{\veps_1-\veps_{h}} ~&:&~
\cc^-_{h-1} ~\lra ~ \cc^-_h ~=~ \cc^-_{h+1} \nn\\
\cd^{m_{1}}_{\veps_1+\veps_{h+1}} ~&:&~
\cc^-_h ~=~ \cc^-_{h+1} ~\lra ~  \cc^+_h ~=~ \cc^+_{h+1} \nn\\
\cd^{m_{3}}_{\veps_1+\veps_{h}} ~&:&~
\cc^+_h ~=~ \cc^+_{h+1} ~\lra ~ \cc^+_{h-1}
\een

\md

\bu ~~$m_{1}=0$ ~equiv~ $n_{2} = -n_{1}$,
\begin{eqnarray}\label{sgnevone} _1\chi^\pm_1 &=& [
 (0,\, m_2,\, \ldots,\, m_\hh)^\pm
\ ; \  \pm (\ha m_2 +  m_{3,h+1} ) \ ]  \ ,  \\
_1\chi^\pm_2 &=& [   (0,\, m_2,\, \ldots,\, m_{\hh-1},\,
m_{h,\hh+1})^\pm
    \ ; \  \pm (\ha m_2 + m_{3,h}  )  \ ]      \nn\\
_1\chi^\pm_3  &=&  [
(0,\, m_2,  \ldots,  m_{\hh-2},  m_{h-1,\hh},  m_{\hh+1})^\pm
\ ; \  \pm
(\ha m_2 + m_{3,h-1} ) \ ]      \nn\\
... \nn\\
_1\chi^\pm_i  &=&  [
(0,\, m_2,  \ldots,  m_{h-i+1},  m_{h-i+2,h-i+3},  m_{h+4-i}, \ldots, \nn\\
&&\qquad\qquad
m_{\hh},\, m_{\hh+1})^\pm \ ; \  \pm
(\ha m_2 + m_{3,h+2-i} ) \ ]      \nn\\
... \nn\\
 _1\chi^\pm_{\hh-1} &=& [  (0,\,
m_2,\, m_{34},\,m_5,\,\ldots,\, m_{\hh},\, m_{\hh+1})^\pm
\ ; \
\pm  (\ha m_2 + m_{3} )  \ ]   \nn\\
 _1\chi^\pm_{\hh} = _1\chi^\mp_{\hh+1} &=& [  (m_{3},\,
m_2+m_{3},\, m_4,\,\ldots,\, m_{\hh},\, m_{\hh+1})^\pm
\ ; \  \pm  \ha m_2  \ ]     \ , \qquad \spa   \nn\end{eqnarray}
Here there are three differential operators involving physically relevant representations,
cf. Figure 11~:
\eqna{difffon} \cd^{m_{3}}_{\veps_1-\veps_{h}} ~&:&~
\cc^-_{h-1} ~\lra ~ \cc^-_h ~=~ \cc^+_{h+1} \nn\\
\cd^{m_{2}}_{\veps_1-\veps_{h+1}} ~&:&~
\cc^-_h ~=~ \cc^+_{h+1} ~\lra ~  \cc^+_h ~=~ \cc^-_{h+1} \nn\\
\cd^{m_{3}}_{\veps_1+\veps_{h}} ~&:&~
\cc^+_h ~=~ \cc^-_{h+1} ~\lra ~ \cc^+_{h-1}
\een

Note that the last two cases: \eqref{sgnevtwo} and \eqref{sgnevone} are conjugate to each other through
the ~$\cm$~ labels ($_1\chi^\pm_{i}$ has the same expressions for ~$c$~ as $_2\chi^\pm_{i}$, but
the ~$\cm$~ labels are conjugate).


For future reference we summarize the physically relevant pairs  reparametrizing
for more natural presentation and introducing uniform notation ~$_r\chi^\pm_k$~:
\begin{eqnarray}\label{sgnevsum} _r\chi^\pm_1  &=& [
 (m_1,\, \ldots,\, m_\hh)^\pm
\ ; \  \pm (\ha m_{12} +  m_{3,h} ) \ ]\ ,
\quad d^+ \geq 2h-1   \ , ~~d^- \leq 1 \ , \nn\\
_r\chi^\pm_2  &=& [   (m_1,\, \ldots,\,
m_{\hh})^\pm
    \ ; \  \pm (\ha m_{12}+ m_{3,h-1}  )  \ ]  \ ,
\quad d^+ \geq 2h-2   \ , ~~d^- \leq 2 \ ,   \nn\\
... \nn\\
_r\chi^\pm_{j}  &=& [ (m_1,\, \ldots,\, m_\hh)^\pm    \ ; \
\pm   (\ha m_{12} +  m_{3,h+1-j} )   \ ] \ , \nn\\
&& \quad d^+ \geq 2h-j   \ , ~~d^- \leq j \ , \quad 1\leq j \leq h-2\ ,  \\
... \nn\\
_r\chi^\pm_{\hh-1}  &=& [  (m_1,\, \ldots,\, m_{\hh})^\pm \ ; \
\pm  \ha m_{12}   \ ]\ ,
\quad d^+ \geq h+1   \ , ~~d^- \leq h-1 \ ,   \nn\\
_r\chi^\pm_{\hh}  &=& [  (m+m_{2},\,
m_{2},\,m_{3},\,\ldots,\, m_{\hh})^\pm \ ; \  \pm  \ha m  \ ]\ ,
\nn\\ &&\qquad\qquad\qquad \qquad d^+ \geq h+\ha   \ , ~~d^- \leq h-\ha \ , \nn\\
 _r\chi^\pm_{\hh+1}  &=& [  (m_{2},\,
m+m_{2},\, m_3,\,\ldots,\, m_{\hh} )^\pm \ ; \  \pm  \ha
m  \ ] \ , \nn\\ &&\qquad\qquad\qquad \qquad
  d^+ \geq h+\ha   \ , ~~d^- \leq h-\ha \ .  \nn\end{eqnarray}

 Note a last reduction obtained by setting ~$m=0$~ when the last two pairs in \eqref{sgnevsum} coincide and become further a singlet
 (being ~$\cm$~ self-conjugate):
\eqn{single}
_r\chi^s   = [  m_{2},\, m_{2},\, m_3,\,\ldots,\, m_{\hh}  \ ; \  0  \ ] \ ,
\quad d =h \ .
\eeq

\subsection{Special cases for ~$p+q$~ even}

The ERS  ~$\chi^-_1$~ are the only ones in the multiplet that contain as irreducible subrepresentations
the finite-dimensional   irreducible representations of ~$\cg$. More precisely, the ER ~$\chi^-_1$~
contains the finite-dimensional   irreducible representation of ~$\cg$~
with signature ~$(m_1, \ldots, m_{h+1})$. (Certainly, the latter is non-unitary
except the case of the trivial one-dimensional  obtained for ~$m_i=1$, ~$\forall\, i$.)

Another important case is the ER with signature
~$\chi^+_1\,$. For ~$pq\in 2\bbn$~ it contains a unitary discrete series representation  of $so(p,q)$ realized on an
invariant subspace ~$\cd$~ of
the ER ~$\chi^+_1\,$. That subspace is annihilated by the KS operator
~$G^-_1\,$,\ and is the image of the KS operator ~$G^+_1\,$.

Furthermore when ~$p>q=2$~ the invariant subspace ~$\cd$~ is the direct sum of two subspaces
~$\cd ~=~ \cd^+ \oplus \cd^-$, in which are realized a
 {\it holomorphic discrete series representation} and its conjugate
  {\it   anti-holomorphic discrete
series representation}, resp.
Note that the corresponding lowest weight GVM is infinitesimally
equivalent only to the holomorphic discrete series, while the
conjugate highest weight GVM is infinitesimally equivalent to the
anti-holomorphic discrete series.

Thus, the signatures of the (holomorphic) discrete series are:
\eqn{holoe}
\chi^+_1 = [  m_1,\, \ldots,\, m_\hh
\ ; \  d ~=~ h +  \ha m_{12}  + m_{3,h} ~+~ \nu \ ]\ , \quad \nu\in\bbn
\eeq

More (non-holomorphic) discrete series representations are contained in
 ~$\chi^+_k$~ for ~$1<k\leq h+1$.

The next important case of positive energy UIRs are the
~{\it limits of (holomorphic) discrete series}~
which are contained in the reduced case \eqref{sgnevsum}:
\eqn{limhoe}
_r\chi^+_1 ~=~  [  m_1,\, \ldots,\,
m_\hh \ ; \ d ~=~ h +  \ha m_{12} + m_{3,h} \ ] \eeq
(with conformal weight obtained from \eqref{holoe} as "limit" for $\nu=0$).

Further we discuss the so called ~{\it first reduction points} (FRP).
These are the boundary values of ~$d$~ from below
of the positive energy UIRs.  Most of the FRPs  are
contained in ~$\chi^+_{h}$, cf. \eqref{sgnev}, which we give with suitable reparametrization:
\eqn{frpe} \chi^+_{h} ~=~ [ m_1,\, m_2,\, \ldots,\, m_{\hh}\ ; \  d ~=~ h  +   \ha m_{12} -1\ ]\ ,
\qquad m_1,m_2 \geq 2 \ .\eeq
Some   FRP cases when only one of ~$m_1,m_2$~ is equal to ~$1$~
are found in  ~$\chi^\pm_{h+1}$~:
\begin{eqnarray}\label{holeone}
\chi^-_{\hh+1} ~&=&~ [  m_1,\, 1,\,m_3,\, \ldots,\, m_{\hh} \ ; ~  d ~=~ h +\ha (m_1-3) \ ]\ ,
\quad m_1\geq 3 \ , \nn\\
\chi^+_{\hh+1} ~&=&~ [ 1,\, m_2,\, \ldots,\, m_{\hh}\ ; ~ d ~=~ h +\ha  (m_2-3)  \ ] \ ,
\quad m_2\geq 3 \ .
\end{eqnarray}
Finally the last three FRP cases ~$(m_1,m_2) = (1,1), (2,1), (1,2)$~ are found in
~$_r\chi^-_{k=,h-1h,h+1}$~:
\begin{eqnarray}\label{holetwoone}
_r\chi^-_{\hh-1} ~&=&~ [  1,\, 1,\,m_3,\, \ldots,\, m_{\hh} \ ; ~  d ~=~ h -1  \ ]\ ,  \nn\\
_r\chi^-_{\hh} ~&=&~ [  2,\, 1,\,m_3,\, \ldots,\, m_{\hh} \ ; ~  d ~=~ h -\ha  \ ]\ ,  \nn\\
_r\chi^-_{\hh+1} ~&=&~ [  1,\, 2,\,m_3,\, \ldots,\, m_{\hh} \ ; ~  d ~=~ h -\ha  \ ]\ .
\end{eqnarray}

\subsection{Minimal irreps for $p+q$ even}

The minimal irreps in this case happen to be related to the ERs in the reduced multiplets.
We define the minimal irreps ~$L_\L$~ as positive energy UIRs which
involve the lowest dimensional representation of ~$\cm$.  Besides the
signature we display the equations that are obeyed by the functions of the irrep.
Typically, the irrep is the intersection of the kernel of the corresponding KS operator
~$G^+$~ and of one or two intertwining differential operators that were already displayed
in the subsection on reduced multiplets.

Below we denote by ~$_rL^\pm_i$~ the irreducible subrepresentation of the ER ~$_r\cc^\pm_i\,$.
The list is:
\eqna{sgnevsumz} _r\chi^-_1  &=& [
 (1,\, \ldots,\, 1)\ ; ~   d ~=~ 1  \ ]\ , \\
&& _rL^-_1 ~=~ \{\ \vf \in\ _r\cc^-_1  ~:~
\cd^{1}_{\veps_1-\veps_{3}}\, \vf ~=~ 0 \ ,
\quad G^+_1\, \vf ~=~ 0 \ \} \ , \nn\\
_r\chi^-_2  &=& [   (1,\, \ldots,\,1)
    \ ; ~  d ~=~ 2   \ ]  \ ,    \nn\\
&& _rL^-_2 ~=~ \{\ \vf \in\ _r\cc^-_2  ~:~
\cd^{1}_{\veps_1-\veps_{4}}\, \vf ~=~ 0 \ ,
\quad G^+_2\, \vf ~=~ 0 \ \} \ , \nn\\
... \nn\\
_r\chi^-_{j}  &=& [ (1,\, \ldots,\,1)    \ ; \
 d ~=~ j    \ ] \ , \quad 1\leq j \leq h-2\ ,\nn\\
&&_rL^-_j ~=~ \{\ \vf \in\ _r\cc^-_j  ~:~
\cd^{1}_{\veps_1-\veps_{j+2}}\, \vf ~=~ 0 \ ,
\quad G^+_j\, \vf ~=~ 0 \ \} \ ,\nn\\
... \nn\\
_r\chi^-_{\hh-1}  &=& [ (1,\, \ldots,\, 1)    \ ; \
d_{\rm FRP} = h-1    \ ] \ ,    \nn\\
&&_rL^-_{h-1} ~=~ \{\ \vf \in\ _r\cc^-_{h-1}  ~:~
\cd^{1}_{\veps_1-\veps_{h+1}}\, \vf ~=~ 0 \ , \quad
\cd^{1}_{\veps_1+\veps_{h+1}}\, \vf ~=~ 0 \ , \nn\\
&&\qquad\qquad\qquad G^+_{h-1}\, \vf ~=~ 0 \ \}  \ ,    \nn\\
_r\chi^-_{\hh}  &=& [  (2,\,1,\, \ldots,\, 1) \ ; \  d^-_{\rm FRP} = h-\ha    \ ]
\ , \nn\\
&&_rL^-_h ~=~ \{\ \vf \in\ _r\cc^-_h  ~:~
\cd^{2}_{\veps_1+\veps_{h+1}}\, \vf ~=~ 0 \   \} , \quad
 G^+_h ~\sim~ \cd^{1}_{\veps_1+\veps_{h+1}}  \ , \nn\\
 _r\chi^-_{\hh+1}  &=& [  (1,\,
2,\, 1,\, \ldots,\, 1) \ ; \  d_{\rm FRP} = h-\ha
  \ ]\ , \nn\\
&&_rL^-_{h+1} ~=~ \{\ \vf \in\ _r\cc^-_{h+1}  ~:~
\cd^{2}_{\veps_1-\veps_{h+1}}\, \vf ~=~ 0 \   \} ,\nn\\ && \qquad
 G^+_{h+1} ~\sim~ \cd^{1}_{\veps_1-\veps_{h+1}}  \ ,
   \nn\end{eqnarray}
where we have indicated (in the last two cases) the degeneration of KS integral
operators to differential operators.

We see in \eqref{sgnevsumz}
that for ~$h\geq 3$~ there are discrete unitary points ~{\it below}~ the FRPs.\footnote{Thus,
the most famous case ~$so(4,2)$~ is excluded.} For fixed ~$h\geq 3$~ these are
in ~$_r\chi^-_j$~ with conformal weight ~$d=j$~ (and trivial $\cm$ inducing irreps)
for ~$j = 1,\ldots, h-2$. Furthermore, as evident from \eqref{sgnevsum}
for ~$h\geq 4$~  there are discrete unitary points below those displayed.
For fixed ~$h\geq 4$~ these are
in ~$_r\chi^-_j$~ with conformal weight ~$1\leq d <j$~ (and non-trivial ~$\cm$~ inducing irreps)
for ~$j = 2,\ldots, h-2$. It seems that all this picture is consistent with \cite{EHW}.
More details  will be given elsewhere.

\bigskip

\nt {\bf Singular vectors needed for the invariant differential operators:}

\md

The necessary  cases are:
\eqna{singev}
&& \veps_1 - \veps_j ~=~ \a_{h+3-j} + \cdots + \a_{h+1}\ , ~~2\leq j \leq h+1 \ ,
\nn\\
&& \veps_1 + \veps_{h+1} ~=~ \a_1 + \a_3 + \cdots + \a_{h+1} \ . \een
These are  roots of ~$sl(n)$~ subalgebras ($n <h+1$). Thus, we can use
f-la \eqref{singone} after suitable change of enumeration.

\section{Figures}

\fig{}{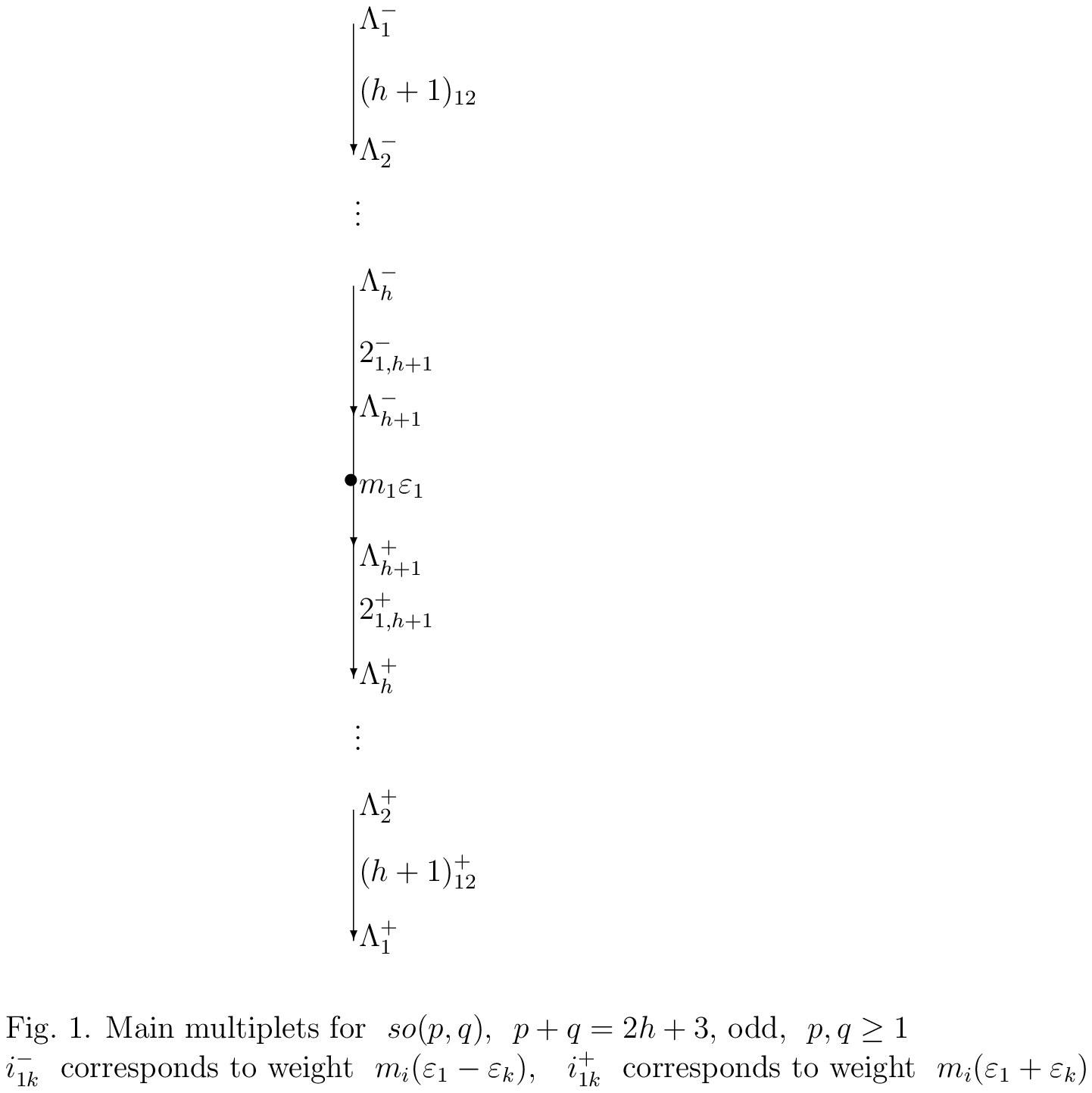}{9cm}
\fig{}{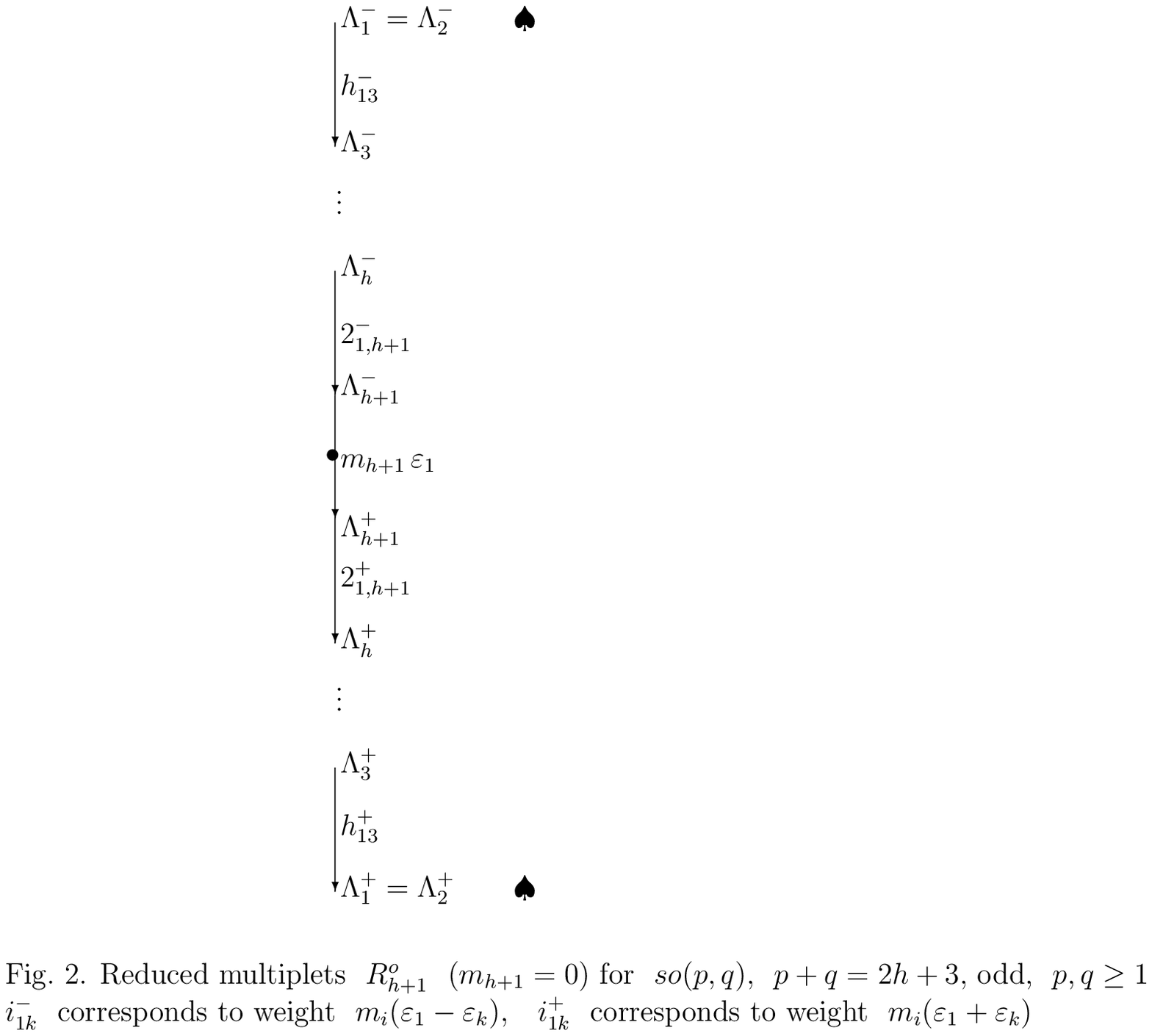}{9cm}
\fig{}{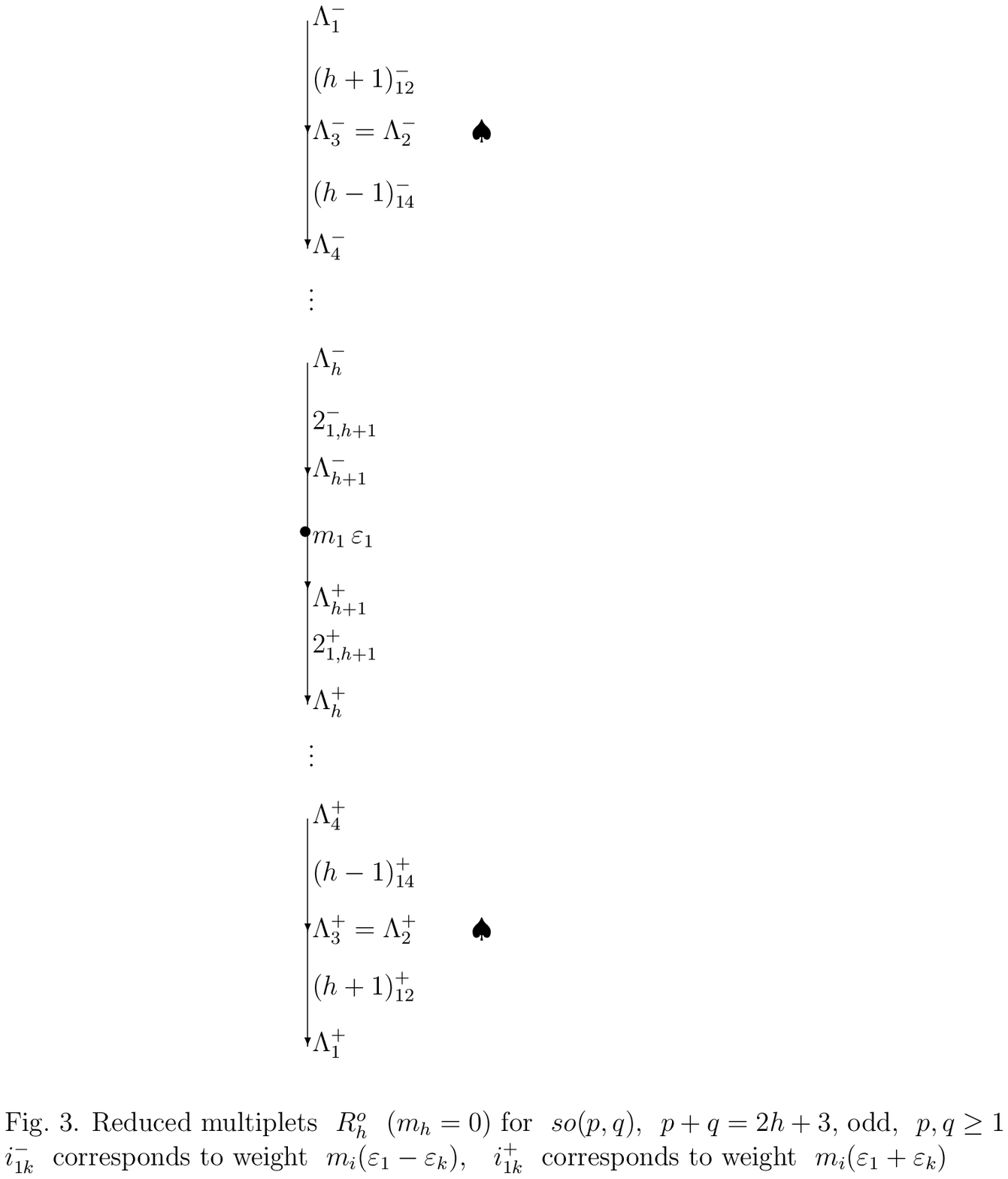}{9cm}
\fig{}{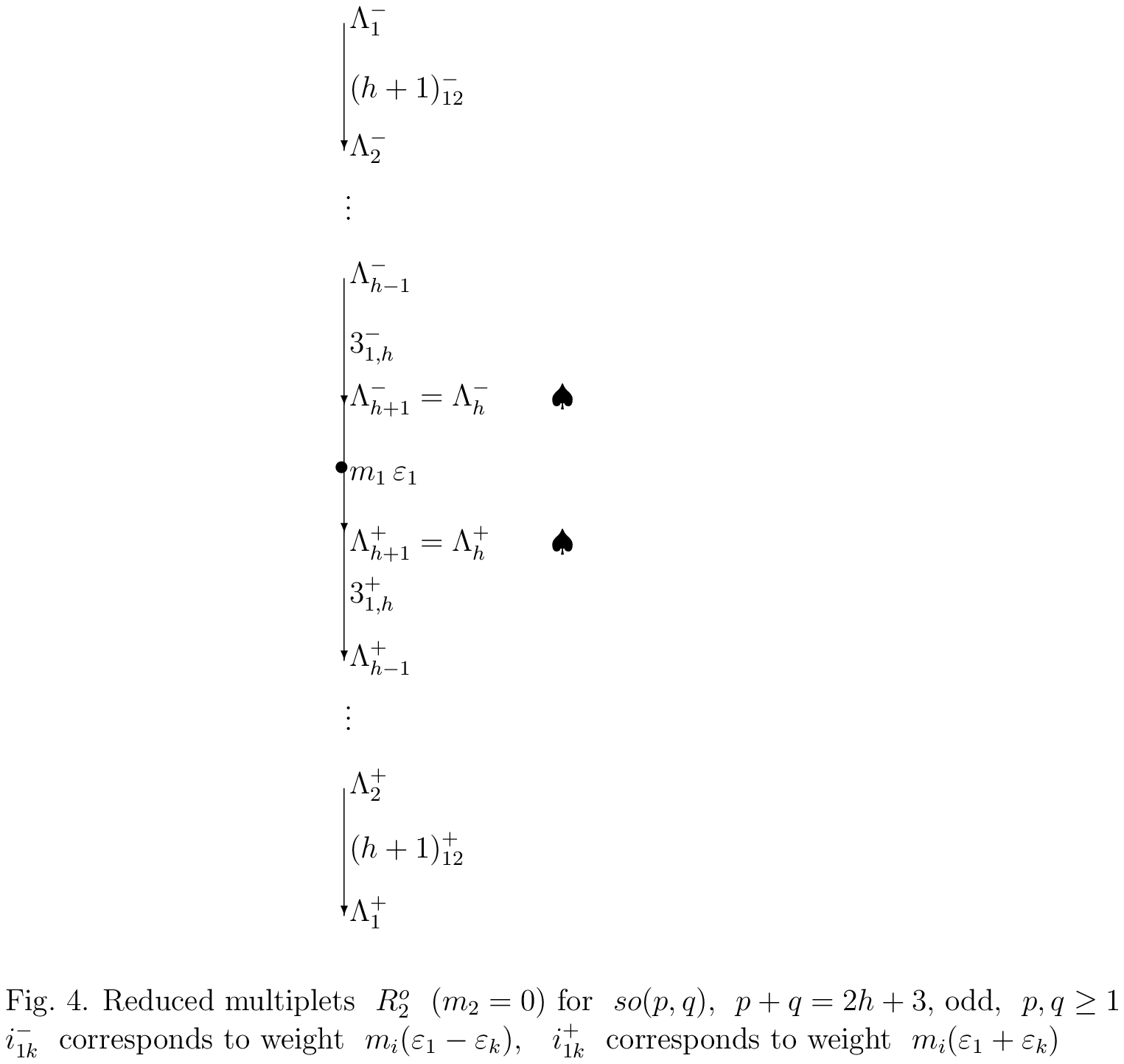}{9cm}
\fig{}{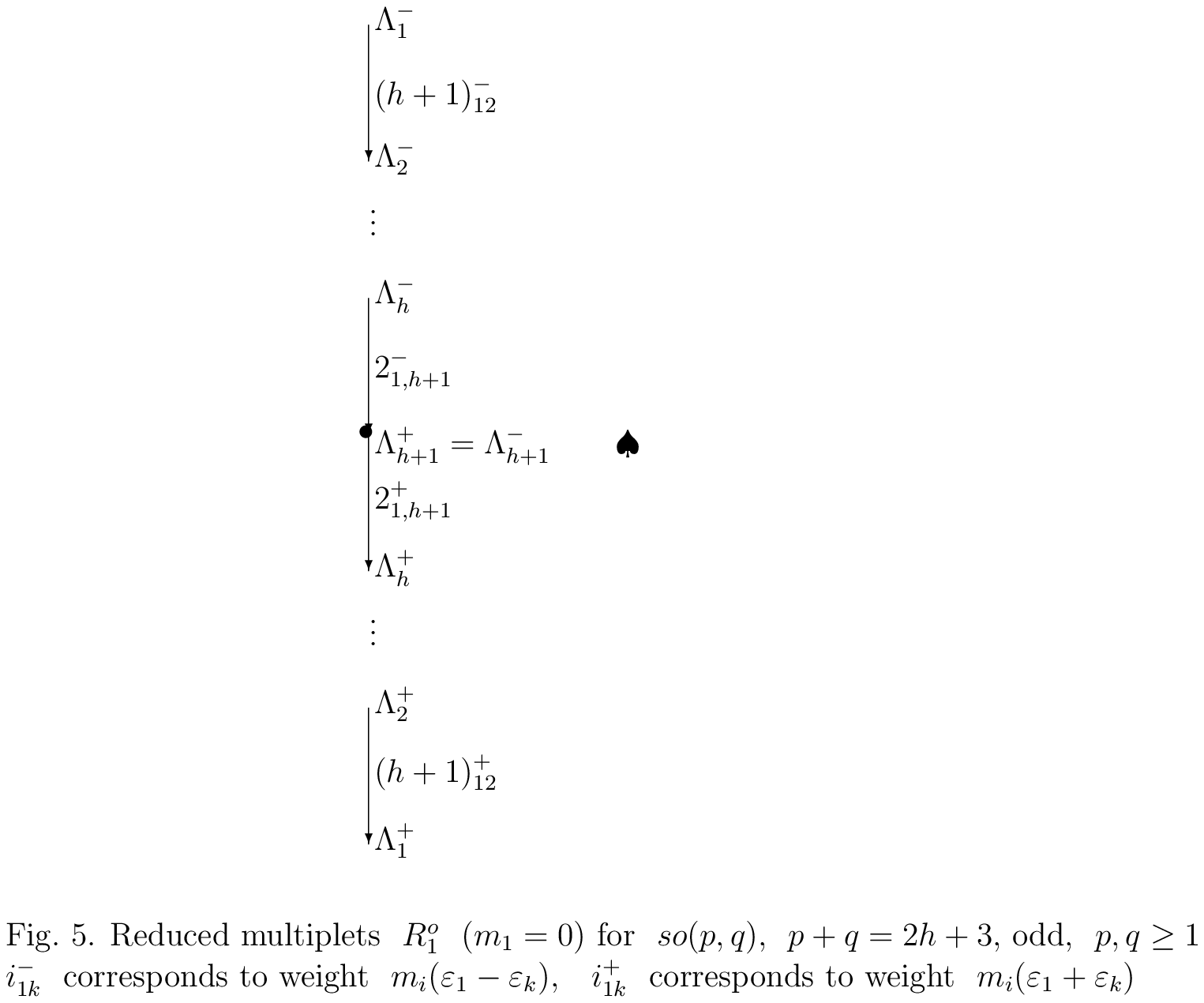}{9cm}

\fig{}{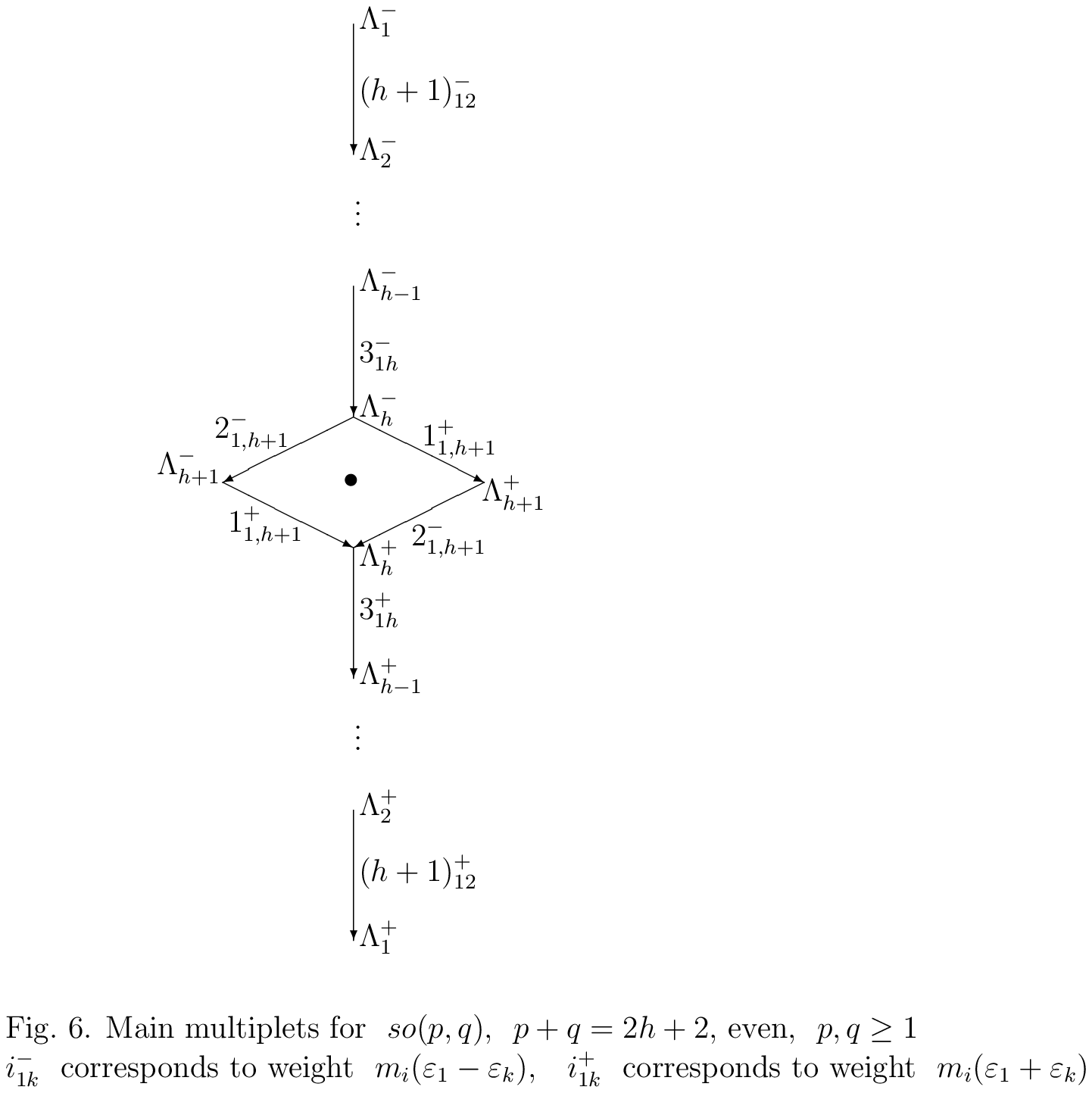}{9cm}
\fig{}{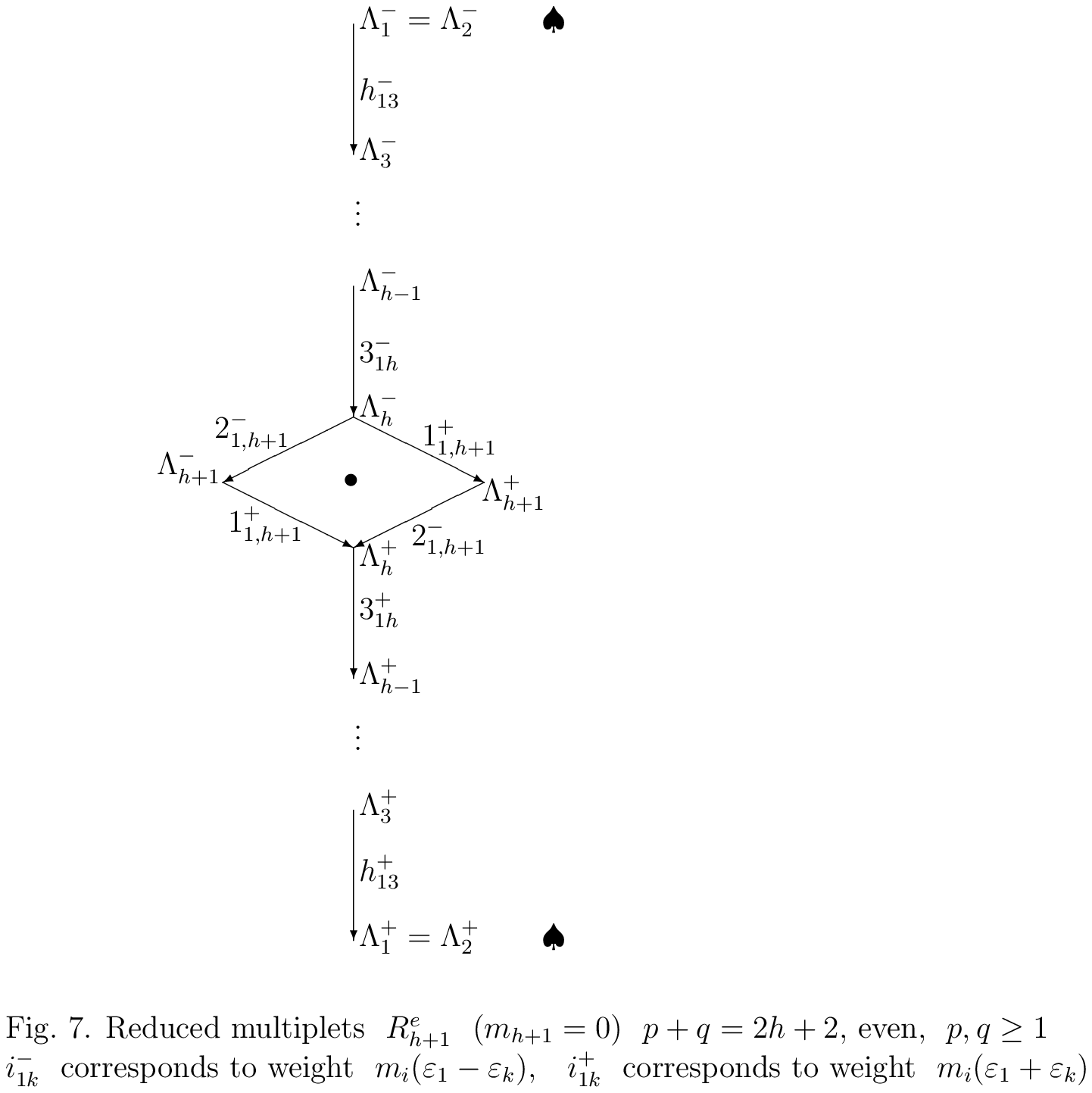}{9cm}
\fig{}{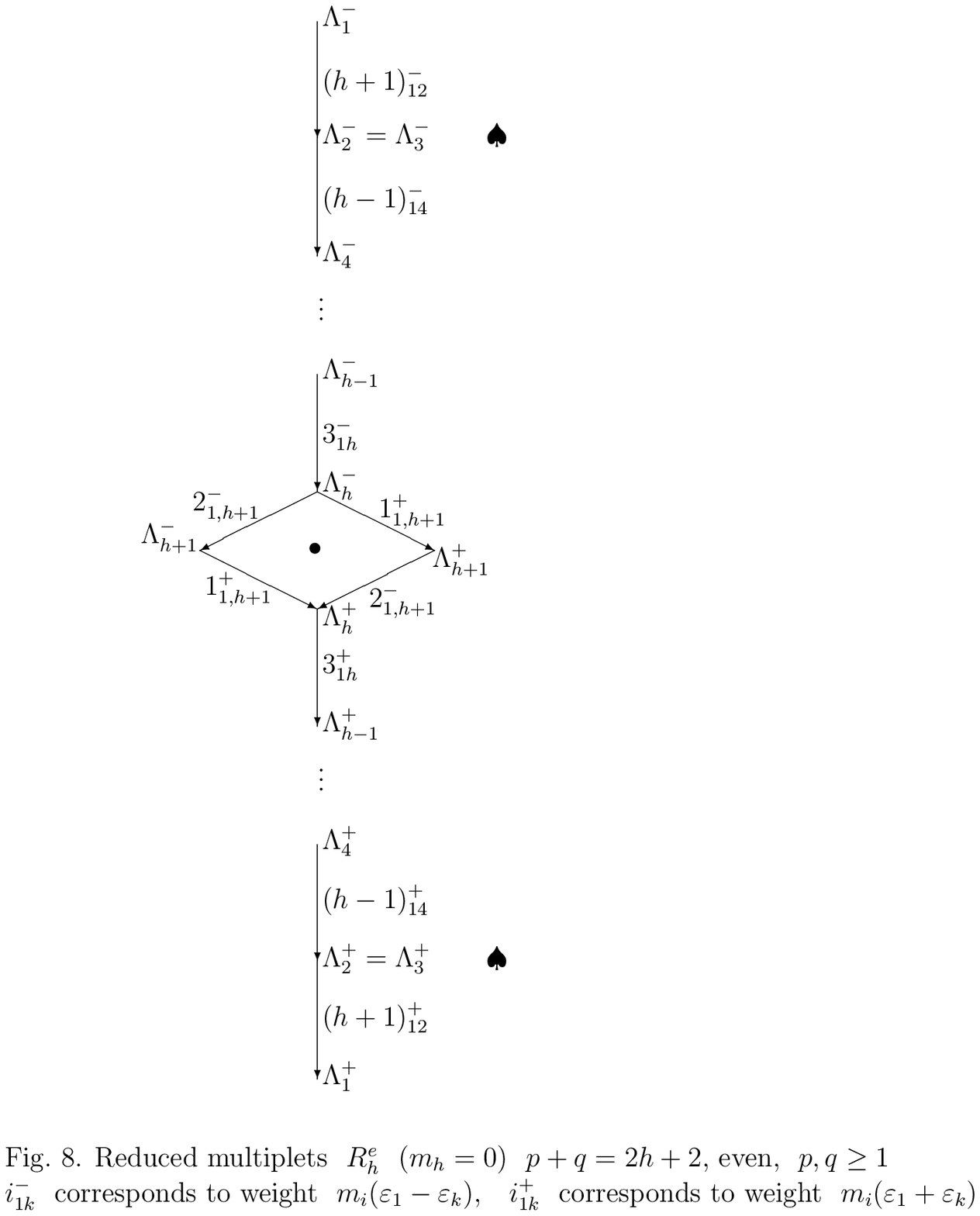}{9cm}
\fig{}{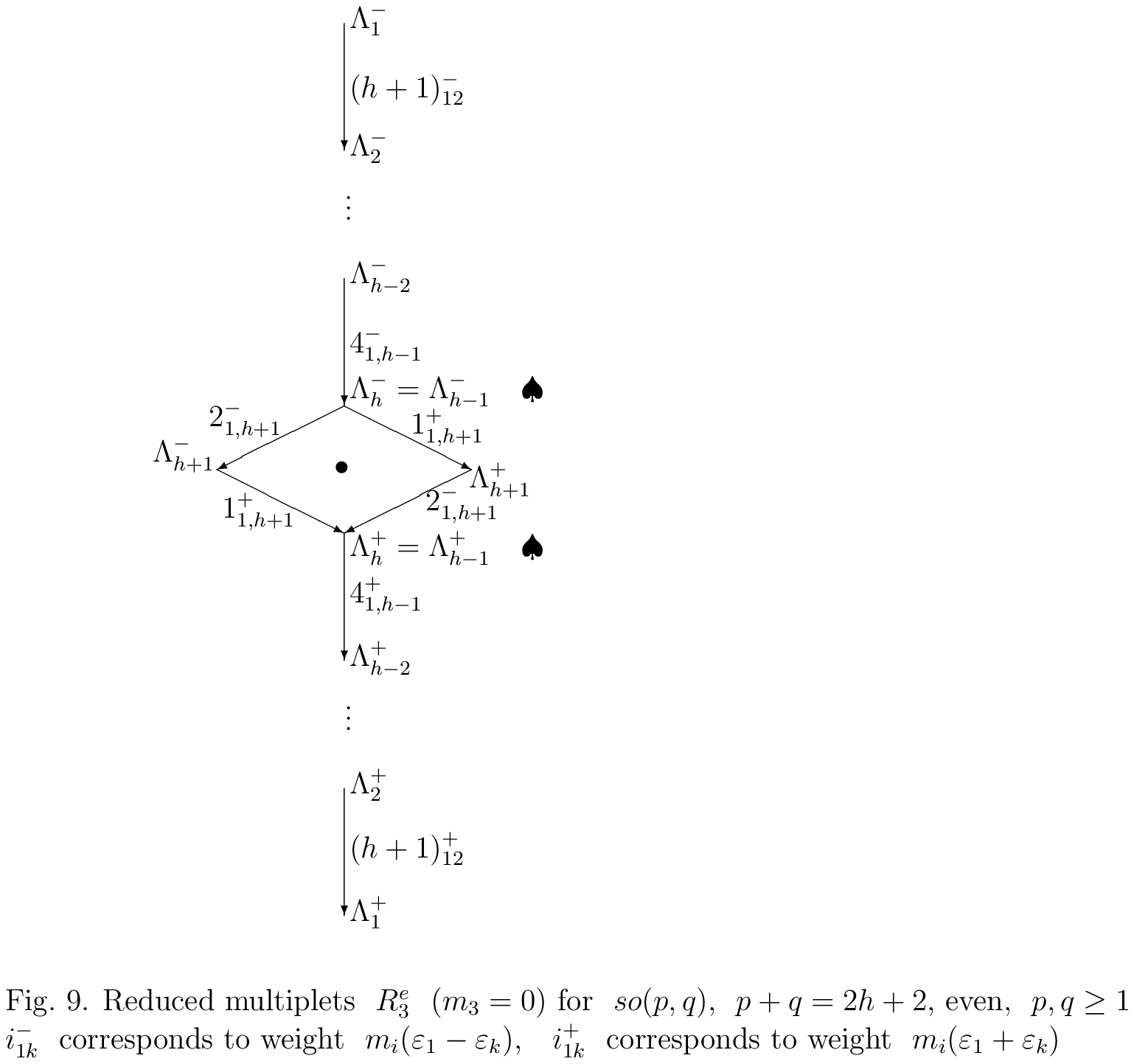}{9cm}
\fig{}{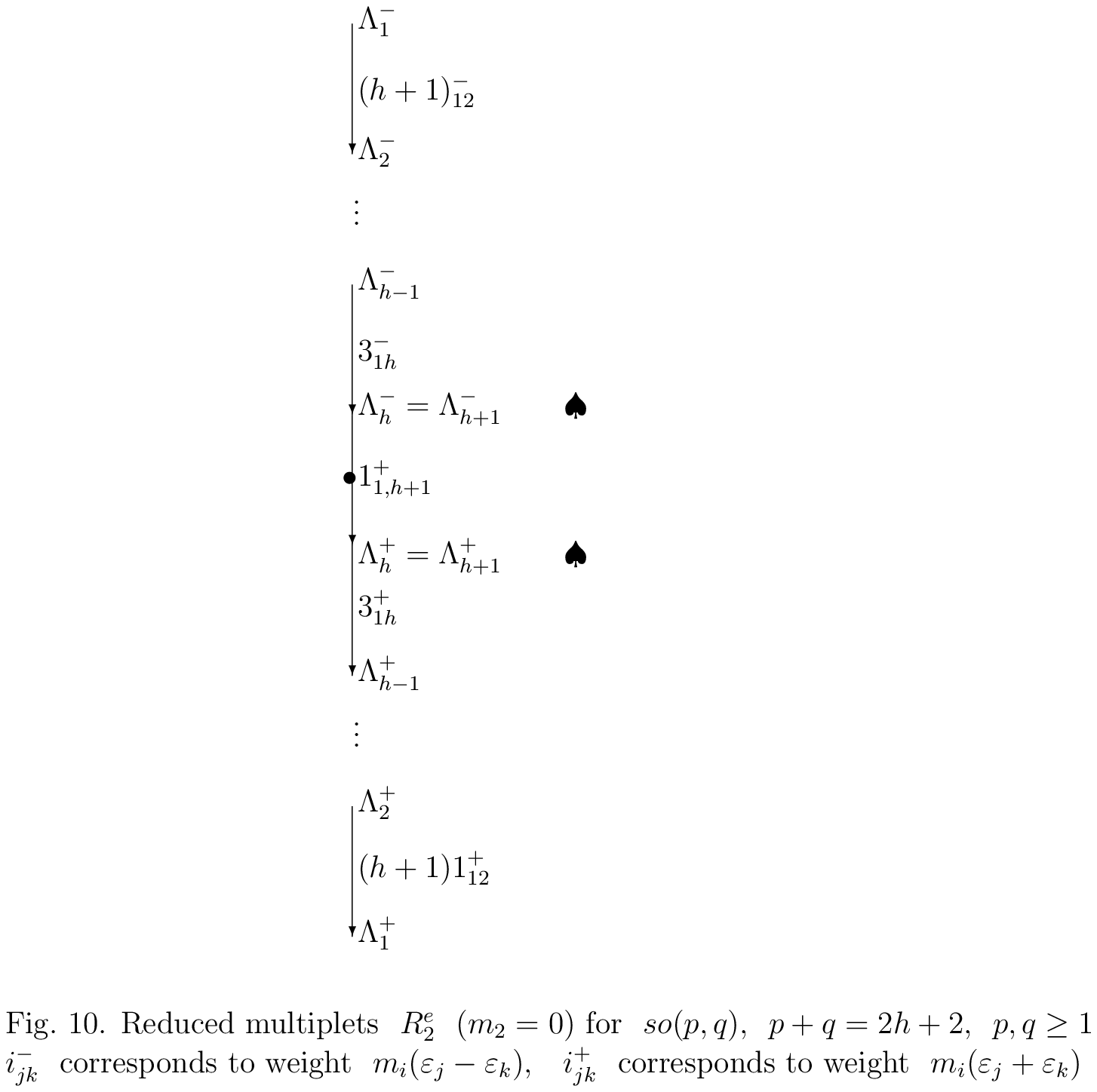}{9cm}
\fig{}{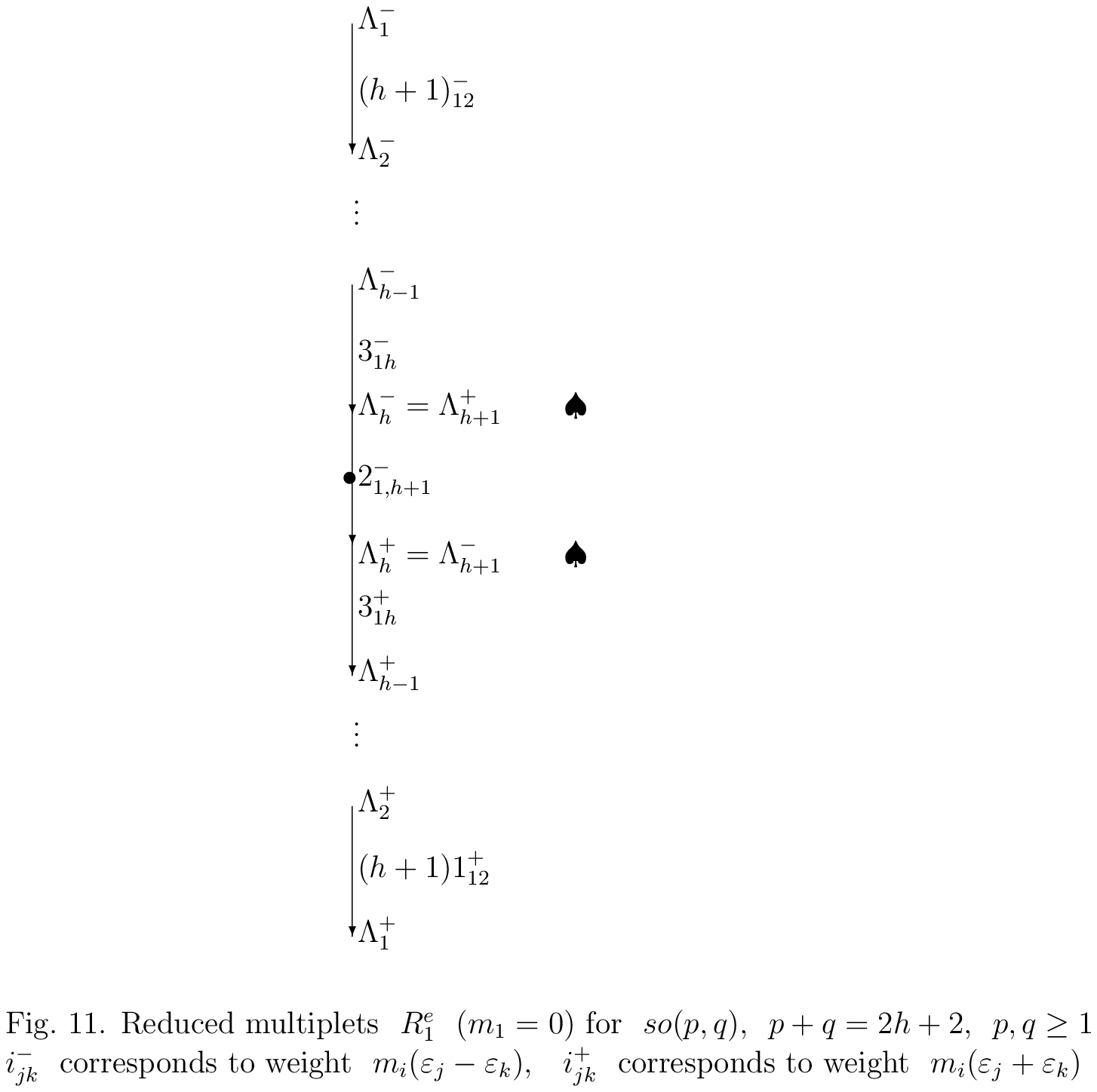}{9cm}

\end{document}